\documentclass[11pt]{amsart}
\usepackage{amssymb}

\newtheorem{theo}{Theorem} 

\newtheorem{lemma}{Lemma}[section]
\newtheorem{prop}[lemma]{Proposition}
\newtheorem{corol}[lemma]{Corollary}
\newtheorem{claim}[lemma]{Claim}

\theoremstyle{remark}
\newtheorem{remark}[lemma]{Remark}
\newtheorem*{remark*}{Remark}

\theoremstyle{definition}

\newcommand{\CC}{\mathbb{C}}

\newcommand{\RR}{\mathbb{R}}
\newcommand{\ZZ}{\mathbb{Z}}

\newcommand{\eps}{\varepsilon}

\newcommand{\LLL}{\mathcal{L}}
\newcommand{\MMM}{\mathcal{M}}

\newcommand{\SSS}{\mathcal{S}}
\newcommand{\UUU}{\mathcal{U}}
\newcommand{\VVV}{\mathcal{V}}

\newcommand{\YYY}{\mathcal{Y}}
\newcommand{\ZZZ}{\mathcal{Z}}

\newcommand{\dd}{\mathsf{d}}

\newcommand{\tth}{\tilde{h}}

\newcommand{\tU}{\widetilde{U}}

\newcommand{\cR}{\mathbb{R}}

\newcommand{\cN}{\mathbb{N}}

\newcommand{\ep}{\epsilon}

\newcommand{\lap}{\triangle}

\newcommand{\tv}{\tilde{v}}

\newcommand{\wtQ}{\widetilde{Q}}

\newcommand{\vbar}{\overline{v}}

\newcommand{\ubar}{\overline{u}}

\newcommand{\fbar}{\overline{f}}

\DeclareMathOperator{\re}{Re}
\DeclareMathOperator{\im}{Im}

\DeclareMathOperator{\ess}{ess}

\DeclareMathOperator{\vect}{span}

\newcommand{\ds}{\displaystyle}

\numberwithin{equation}{section} 

\title[3D cubic NLS]{Threshold solutions for the focusing 3d cubic Schr\"odinger
equation}
\author{Thomas Duyckaerts}
\email{thomas.duyckaerts@u-cergy.fr}
\address{Universit\'e de Cergy-Pontoise, UMR CNRS 8088, France}
\author{Svetlana Roudenko}
\email{svetlana@math.asu.edu}
\address{Arizona State University}
\begin{document}

\maketitle

\begin{abstract}
We study the focusing 3d cubic NLS equation with $H^1$ data at the
mass-energy threshold, namely, when $M[u_0]E[u_0]=M[Q]E[Q]$. In
\cite{HoRo07a}, \cite{HoRo07b} and \cite{DuHoRo07P}, the behavior of
solutions (i.e., scattering and blow up in finite time) is
classified when $M[u_0]E[u_0]<M[Q]E[Q]$. In this paper, we first
exhibit 3 special solutions: $e^{it} Q$ and $Q^\pm$, here $Q$ is the
ground state, and $Q^\pm$ exponentially approach the ground state
solution in the positive time direction and $Q^+$ having finite time
blow up and $Q^-$ scattering in the negative time direction.
Secondly, we classify solutions at this threshold and obtain that up
to $\dot{H}^{1/2}$ symmetries, they behave exactly as the above
three special solutions, or scatter and blow up in both time
directions as the solutions below the mass-energy threshold. These
results are obtained by studying the spectral properties of the
linearized Schr\"odinger operator in this mass-supercritical case,
establishing relevant modulational stability and careful analysis of
the exponentially decaying solutions to the linearized equation.
\end{abstract}

\section{Introduction}

We consider the 3d focusing cubic nonlinear Schr\"odinger (NLS)
equation on a time interval $I \subset \cR$ ($0 \in I$)
\begin{equation}
 \label{CP}
\left\{
\begin{gathered}
i\partial_t u +\Delta u+|u|^{2}u=0,\quad (x,t)\in \RR^3\times I\\
u_{\restriction t=0}=u_0\in H^1(\RR^3).
\end{gathered}\right.
\end{equation}
The Cauchy problem \eqref{CP} is locally wellposed in $H^1$, see
\cite{GiVe79}. We denote the forward lifespan by $[0, T_+)$ and the
backward by $(T_-, 0\,]$. If $T_+(u) < + \infty$ [or $T_-(u)
> -\infty$], then $\lim_{t \to T_+} \|u(t)\|_{H^1}=+\infty$,
[respectively, $\lim_{t \to T_-} \|u(t)\|_{H^1}=+\infty$], and it is
said that the solution blows up in finite time.

The solutions of \eqref{CP} satisfy mass, energy and momentum
conservation laws
\begin{align*}
E[u](t)&=\frac{1}{2} \int |\nabla u(x,t)|^2\, dx - \frac{1}{4} \int
|u(x,t)|^{4} \, dx = E[u](0),\\
M[u](t)&=\int |u(x,t)|^2 dx=M[u](0),\\
P[u](t)&=\im \int \overline{u}(x,t)\nabla u(x,t)dx= P[u](0).
\end{align*}
Furthermore, this NLS equation enjoys several invariances. If
$u(t,x)$ is a solution, then
\begin{itemize}
\item[{\bf -}]
by {\it scaling invariance}: so is $\lambda \, u(\lambda x,
\lambda^2 t)$, $\lambda>0$;
\item[{\bf -}]
by {\it spatial translation}: so is $u(x+x_0, t)$ for $x_0 \in
\cR^3$;
\item[{\bf -}]
by {\it time translation}: so is $u(x,t+t_0)$ for $t_0\in \RR$;
\item[{\bf -}]
by {\it phase rotation invariance}: so is $e^{i\theta_0} u$,
$\theta_0 \in \cR$;
\item[{\bf -}]
by {\it time reversal symmetry}: so is $\overline{u(x, -t)}$.
\end{itemize}
Observe that all these transformations leave the
$\dot{H}^{1/2}$-norm and the momentum invariant.
In what follows the solutions will be considered {\it up to the
($\dot{H}^{1/2}$-) symmetries} of this NLS equation meaning up to the above mentioned
invariances.

A transformation of solutions to \eqref{CP}, which does not leave
the $\dot{H}^{1/2}$-norm nor the momentum invariant, is the {\it
Galilean transformation}: if $u$ is a solution, so is
$$
e^{ix\xi_0}e^{-it|\xi_0|^2}u\left(x-2\xi_0t,t\right), \; \xi_0 \in
\cR^3.
$$

Consider a general focusing NLS equation
$$
i \partial_t u +\Delta u + |u|^{p-1}u=0, \quad (x,t) \in \cR^d
\times \cR,
$$
with the nonlinearity $p>1$ and the dimension $d$ such that $0 \leq
s_c \leq 1$, where $s_c = \frac{d}2-\frac2{p-1}$. The case when $s_c
= 0$ is referred to as mass (or $L^2$)-critical, the case when $s_c
= 1$ is called energy-critical, and in our case, the NLS equation in
\eqref{CP} has $s_c = 1/2$, and thus, is referred to as
$\dot{H}^{1/2}$-critical.

The focusing mass-critical NLS equation (for example, cubic NLS in
2d) with $H^1$ initial data was originally studied by Weinstein
\cite{We83}, who showed that
there exists a sharp threshold,
which splits the behavior of solutions: (i) if
$M[u] < M[Q]$, then the solution exists globally in time and (ii)
if $M[u] \geq M[Q]$, then the solution may blow up in finite time. Here,
$Q$ is the {\it ground state} solution of $-Q+\Delta Q+|Q|^{4/d}
Q=0$, $Q=Q(r)$, $r=|x|$, $x \in \cR^d$.
In the first case the scattering was known for the initial
data in $L^2$ and of finite variance (it follows from the
pseudoconformal conservation law, e.g., see \cite{We83},
\cite{Bo99JA}). The scattering for radially symmetric solutions with
$L^2$ initial data was recently established in \cite{KTV} for 2d,
and in \cite{KVZ} for higher dimensions. For general $L^2$ initial
data scattering is still an open question.

Note that the solution $u(x,t) = e^{it} Q(x)$ (it has $M[u] = M[Q]$)
exists globally in time (in fact, it is time-periodic), but does not
scatter. Under the pseudoconformal transformation this solution can
be mapped into a finite time blow up solution (of the same mass).
Merle has shown that all finite time blow up solutions of minimal
mass $M[u] = M[Q]$ are pseudoconformal images (up to phase,
translation, scaling and Galilean invariances) of $e^{it}Q(x)$, see
\cite{Me92} for radial $H^1$ data with finite variance and
\cite{Me93} for general $H^1$ data. Furthermore, he characterized
all $H^1$ solutions of finite variance with the threshold mass
$M[u]=M[Q]$: a solution can be a scaled version of the time periodic
solution $e^{it}Q(x)$, or a blow up solution which is a
pseudoconformal image of $e^{it}Q(x)$ (a ``self-similar solution"),
or a globally defined solution with quadratically decaying in time
$L^{4/d+2}$ norm which implies scattering as $t \to \pm\infty$.

The focusing energy critical NLS equation (for example, cubic NLS in
4d or quintic NLS in 3d) was recently studied by Kenig-Merle
\cite{KeMe06}. They showed that (in dimensions $3, 4$ and $5$) a
sharp splitting takes place for the Cauchy problem with
$\dot{H}^1_{rad}$ initial data and an {\it a priori} condition
$E[u_0] < E[W]$: (i) if $\|\nabla u_0 \|_{L^2} < \|\nabla W
\|_{L^2}$, then the solution exists globally in time, moreover, it
scatters; (ii) if $\| \nabla u_0 \|_{L^2} > \|\nabla W \|_{L^2}$,
and $u_0 \in L^2$, then finite time blow up occurs\footnote{Under
the above a priori condition the gradients of $u_0$ and $W$ can not
be equal.}. Here, $W$ is the stationary solution of \eqref{CP} in
$\dot{H}^1$, given explicitly by $W(r)= 1/\big(1+ \frac{r^2}{d(d -2)
}\big)^{(d -2)/2}$, $r=|x|$, $x \in \cR^d$. A similar result (but
not necessarily for radial initial data) is established by the same
authors for the energy-critical focusing nonlinear wave (NLW)
equation in \cite{KeMe06b}.

Observe that the above characterization is obtained only if $E[u_0]
< E[W]$. What happens if this condition is removed? The case of the
critical level of energy, i.e., $E[u_0] = E[W]$,
was recently studied by Duyckaerts-Merle in \cite{DuMe07}. Richer
dynamics for the behavior of solutions as $t \to \pm \infty$ are
exhibited. Besides the stationary solution $W$ which exists globally
but does not scatter, there are two more special solutions $W^-$ and
$W^+$ which approach $W$ in $\dot{H}^1$ in one time direction, but
in the opposite time direction $W^-$ scatters and $W^+$ blows up in
finite time\footnote{The blow up is shown for $d=5$ and conjectured
for other dimensions.}. The deciding
factor is the gradient size:
$\|\nabla W^- \|_{L^2} < \|\nabla W \|_{L^2}$ and $\|\nabla W^+
\|_{L^2} > \|\nabla W \|_{L^2}$. Moreover, the classification of all
(radial) solutions at the energy critical level is given (up to the
symmetries of the equation): if $\|\nabla u_0 \|_{L^2} < \|\nabla W
\|_{L^2}$, then the solution is $W^-$ , if $\|\nabla u_0 \|_{L^2} >
\|\nabla W \|_{L^2}$ (with the additional technical assumption that
$u_0$ is in $L^2$), then the solution is $W^+$, and when $\|\nabla
u_0 \|_{L^2} = \|\nabla W \|_{L^2}$, then the solution is $W$. A
similar result is obtained for the energy-critical focusing NLW
equation for general initial data in Duyckaerts-Merle
\cite{DuMe07Pb}.

The results on global existence or finite time blow up for the mass
critical NLS and energy critical NLS equations can be linked by
studying the NLS equation with $0<s_c<1$ considered in
Holmer-Roudenko \cite[Section 2]{HoRo07a}, see also \cite[Section
7]{HoRo07b}. For the purpose of this paper we will state the result
only for \eqref{CP}, i.e., when $s_c = 1/2$.  The scattering result
in the following theorem was established initially for the radial
$H^1$ data in \cite{HoRo07b} and the radiality assumption was
removed in \cite{DuHoRo07P}.

Let $Q$ be the ground state, that is the unique positive radial
solution of the equation $-Q+\Delta Q+|Q|^2Q=0$ (see Subsection
\ref{SS:Q} for the details).
\begin{theo}[\cite{HoRo07a, HoRo07b, DuHoRo07P}]
 \label{T:BlowUp}
Let $u$ be an $H^1$ solution to \eqref{CP}. Suppose
\begin{equation}
 \label{E:MEless}
M[u_0]E[u_0] < M[Q]E[Q].
\end{equation}

If $\|u_0\|_{L^2}\|\nabla u_0\|_{L^2}<\|Q\|_{L^2}\|\nabla
Q\|_{L^2}$, then $\|u_0\|_{L^2}\|\nabla u (t) \|_{L^2} < \|Q\|_{L^2}
\|\nabla Q\|_{L^2}$, and thus, the solution $u$ is globally 
defined; moreover, it scatters in $H^1$.

If $\|u_0\|_{L^2}\|\nabla u_0\|_{L^2} > \|Q\|_{L^2}\|\nabla
Q\|_{L^2}$, then $\|u_0\|_{L^2} \|\nabla u(t)\|_{L^2} > \|Q\|_{L^2}
\|\nabla Q\|_{L^2}$ and, if either $u_0$ is radial or has a finite
variance, i.e., $|x|u_0 \in L^2$, then the solution $u$ blows up in
finite time.
\end{theo}
As in the preceding cases, the determining quantities  $M[u]E[u]$
and $\|u_0\|_{L^2}\|\nabla u_0\|_{L^2}$ are invariant by the scaling
of the equation.

Techniques employed above are based on the approach of Kenig-Merle
in \cite{KeMe06} and \cite{KeMe06b}. In particular, scattering is
established by the profile decomposition method by G\'erard
\cite{Ge98}, which is a refinement of the concentration-compactness
method of P.-L. Lions \cite{Li84IHPb,Li85Reb}; see previous
applications to NLS by Merle-Vega \cite{MeVe98} and Keraani
\cite{Ke01}, and to NLW equation by Bahouri-G\'erard \cite{BaGe97}.
For other recent applications of profile decomposition we refer the
reader to the works of G\'erard \cite{Ge08P} on the 3d cubic wave
equation, and of Kenig-Merle \cite{KeMe07P}, who established
scattering of solutions for the defocusing cubic NLS in 3d (equation
\eqref{CP} with a minus sign in front of the nonlinearity) with
$\dot{H}^{1/2}$ initial data provided $\dot{H}^{1/2}$ norm stays
bounded globally in time (see also
Colliander-Keel-Staffilani-Takaoka-Tao \cite{CoKeStTaTa04}
 and references therein for previous results on
the defocusing problem).

Coming back to Theorem \ref{T:BlowUp}, we would like to describe the
behavior of solutions to \eqref{CP} at the ``critical" mass-energy
threshold, i.e., when
\begin{equation}
 \label{E:ME-Q}
M[u]E[u] = M[Q]E[Q].
\end{equation}
First, we establish the existence of special solutions (besides
$e^{it}Q$) at the critical mass-energy threshold.
\begin{theo}
 \label{T:existence}
There exist two radial solutions $Q^+$ and $Q^-$ of \eqref{CP} with initial
conditions $Q^{\pm}_0$ such that $\ds Q^{\pm}_0\in \cap_{s \in \RR}
H^s(\RR^3)$ and
\begin{enumerate}
\item \label{Q+Q-}
$\ds M[Q^+]=M[Q^-]=M[Q],\; E[Q^+]=E[Q^-]=E[Q]$, $[0,+\infty)$ is in
the (time) domain of definition of $Q^{\pm}$ and there exists $e_0>0$ such that
$$
\forall t\geq 0,\quad \left\|Q^{\pm}(t)-e^{it}Q\right\|_{H^1}\leq
Ce^{-e_0t},
$$
\item
$\ds \|\nabla Q^-_0\|_2<\|\nabla Q\|_2$, $Q^-$ is globally defined
and scatters for negative time,
\item
$\ds \|\nabla Q^+_0\|_2>\|\nabla Q\|_2$, and the negative time of existence of $Q^+$ is
finite.
\end{enumerate}
\end{theo}
\begin{remark*}
The best constant $-e_0$ in \eqref{Q+Q-} is given by the negative
eigenvalue of the linearized operator associated to \eqref{CP}
around the periodic solution $e^{it}Q$. Furthermore, the
construction of $Q^{\pm}$ gives an asymptotic expansion in all
Sobolev spaces for all orders of $e^{-e_0t}$ of $Q^{\pm}$ as
$t\rightarrow+\infty$. Such precise information is not available for
negative times. In particular, we are not able to describe the
behavior of $Q^+$ near the blow-up time, except for what is already
known for general blow-up solutions of \eqref{CP} (see
\cite{MeRa06P}, \cite{HoRo07a}).
\end{remark*}
Next, we characterize all solutions at the critical mass-energy
level as follows:

\begin{theo}
\label{T:classification}
Let $u$ be a solution of \eqref{CP} satisfying \eqref{E:ME-Q}.
\begin{enumerate}
\item \label{case_sub}
If $\ds \|\nabla u_0\|_2\|u_0\|_2<\|\nabla Q\|_2\|Q\|_2$, then
either $u$ scatters or $u=Q^-$ up to the symmetries.
\item
 \label{case_critical}
If $\ds \|\nabla u_0\|_2\|u_0\|_2=\|\nabla Q\|_2\|Q\|_2$, then
$u=e^{it}Q$ up to the symmetries.
\item
 \label{case_super}
If $\ds \|\nabla u_0\|_2\|u_0\|_2>\|\nabla Q\|_2\|Q\|_2$ and $u_0$
is radial or of finite variance, then either the interval of
existence of $u$ is of finite length or $u=Q^+$ up to the
symmetries.
\end{enumerate}
\end{theo}
Note that as a consequence of \eqref{E:ME-Q}, the assumptions
$\|\nabla u(t_0)\|_2\|u(t_0)\|_2<\|\nabla Q\|_2\|Q\|_2$ and
$\|\nabla u(t_0)\|_2\|u(t_0)\|_2>\|\nabla Q\|_2\|Q\|_2$ do not
depend on the choice of the initial condition (see \S
\ref{SS:gradient_separation}).
\begin{remark*}
It is worth linking the $\dot{H}^{1/2}$-critical equation \eqref{CP}
with the corresponding mass-critical and energy-critical once again.
The proofs of Theorems \ref{T:existence} and \ref{T:classification}
will show that the behavior of the solutions of \eqref{CP} at the
threshold is very close to the one of the energy critical equation
described in the radial case in \cite{DuMe07}. In particular, in
both cases, the existence of the special solutions $Q^{\pm}$
($W^{\pm}$ in the energy-critical case) derive from the existence of
two real nonzero eigenvalues for the linearized operator around the
periodic solution $e^{it}Q$ (respectively around the stationary
solution $W$). On the other hand, in the mass-critical case, the
only eigenvalue of the linearized operator is $0$ (see \cite{We85}),
and the blow-up solution at the threshold is given by the
pseudo-conformal transformation, which is specific to the
mass-critical equation.
\end{remark*}

We next give a formulation of Theorems \ref{T:BlowUp} and
\ref{T:classification} that takes the Galilean transformation into
account. Let $u$ be a solution of \eqref{CP}. Applying to $u$, as in
\cite[Section 4]{DuHoRo07P}, the Galilean transformation with
parameter $\xi_0=-P[u]/M[u]$, we get a solution $v$ of \eqref{CP}
with zero momentum which is the minimal energy solution among all
Galilean transformations of $u$. Precisely,
$$
M[v] = M[u],\quad E[v]=E[u]-\frac 12\frac{P[u]^2}{M[u]},\quad
\|\nabla v_0\|^2_{L^2}=\|\nabla u_0\|^2_{L^2} -
\frac{P[u_0]^2}{M[u_0]}.
$$
Applying Theorems \ref{T:BlowUp} and \ref{T:classification} to $v$, we get
\begin{theo}
\label{T:classification'}
Let $u$ be a solution of \eqref{CP} satisfying
$$
E[u]M[u]-\frac 12 P[u]^2\leq E[Q]M[Q].
$$
Then
\begin{enumerate}
\item
If $\ds \|\nabla u_0\|_2^2\|u_0\|_2^2-P[u]^2<\|\nabla Q\|_2^2\|Q\|_2^2$, then
either $u$ scatters or $u=Q^-$ up to the symmetries.
\item
If $\ds \|\nabla u_0\|_2^2\|u_0\|_2^2-P[u]^2=\|\nabla Q\|_2^2\|Q\|_2^2$, then
$u=e^{it}Q$ up to the symmetries.
\item \label{CaseBlowUp}
If $\ds \|\nabla u_0\|_2^2\|u_0\|_2^2-P[u]^2>\|\nabla
Q\|_2^2\|Q\|_2^2$ and $u_0$ is radial or of finite variance, then
either the interval of existence of $u$ is of finite length or
$u=Q^+$ up to the symmetries.
\end{enumerate}
\end{theo}
In the preceding theorem, up to the symmetries of the equation means
up to the $\dot{H}^{1/2}$-symmetries \textbf{and} the Galilean
transformation. If $E[u]M[u]\leq E[Q]M[Q]$, our results actually
show that the condition $\|\nabla u_0\|_2^2\|u_0\|_2^2>\|\nabla
Q\|_2^2\|Q\|_2^2$ implies the stronger bound $ \|\nabla
u_0\|_2^2\|u_0\|_2^2-P[u_0]^2>\|\nabla Q\|_2^2\|Q\|_2^2$.
\bigskip

The paper is organized as follows. In the next section we recall the
properties of the ground state $Q$, small data theory for Cauchy
problem \eqref{CP} and the spectral properties of the linearized
(around the ground state solution $e^{it}Q$) Schr\"odinger operator.
Under the condition \eqref{E:ME-Q} we identify a quadratic form
associated to the linearized Schr\"odinger operator which can
measure closeness to $e^{it}Q$ and find subspaces of $H^1$ where
this form is positive, avoiding thus vanishing and negative
directions. In Section \ref{S:spectralsolution} we construct a
family of approximate solutions using the knowledge about the
discrete spectrum of the linearized operator and then with a fixed
point argument produce candidates for the special solutions $Q^-$
and $Q^+$, thus, proving Theorem \ref{T:existence} except for the
negative time behavior of $Q^\pm$. In Section \ref{S:modulation} we
discuss the modulational stability near the ground state solution.
Here, we identify the spatial and phase parameters which control the
variations from $e^{it}Q$ (on the subsets where the above mentioned
quadratic form is positive) while the entire variation being small
in $H^1$ norm. In Sections \ref{S:supercritical} and
\ref{S:subcritical} we study solutions with initial data from
Theorem \ref{T:classification} part \eqref{case_sub} and
\eqref{case_super}, respectively. Our main goal is to obtain
exponential convergence for large (positive) time of the gradient
variation \eqref{E:CVexp} which then will imply exponential
convergence in (positive) time to $e^{it}Q$ (up to spatial
translation and phase rotation), see Lemma \ref{L:expo}. We also
finish Theorem \ref{T:existence} for negative time behavior. In
Section \ref{S:uniqueness} we first analyze exponentially small
solutions of the linearized Schr\"odinger equation and then
establish the uniqueness of special solutions. We finish the section
with the classification of solutions result.
Appendix contains the proof of coercivity of the quadratic form
introduced in Section \ref{S:spectralsolution} where we follow
Weinstein \cite{We85} and a useful inequality, the original idea of
which is due to Banica \cite{Ba04}.

\subsection*{Notation}
Let $\SSS$ denote the space of Schwartz functions, i.e, the
topological space of functions $v$ satisfying
\begin{equation*}
\forall \, \alpha,N,\quad \|v\|_{\alpha,N}:=\sup_{x\in\RR^3}
\Big|(1+|x|)^N\partial_x^{\alpha} v(x)\Big|<\infty,
\end{equation*}
with the topology given by the family of semi-norms
$\|\cdot\|_{\alpha,N}$.

By $H^s$ we denote the usual Sobolev space of smoothness $s$ in
spatial (on $\cR^3$) variable. Let $H^{\infty}= \bigcap_{s \in \cR}
H^s$. We denote by $\|\cdot \|_p$ the $L^p$-norm in spatial
variable.

A couple $(q,r)$ is $\dot{H}^s$-admissible when
$\frac{2}{q}+\frac{3}{r}=\frac{3}{2}-s$ and $2 \leq q,r \leq
\infty$. Consider the following Strichartz norm for functions $u$ of
space and time
$$
\|u\|_{S \left(\dot{H}^{1/2}\right)} = \sup_{\substack{(q,r)\,
\dot{H}^{1/2} - {\rm admissible}\\4^+\leq q\leq \infty,\,3 \leq r
\leq 6^{-}}} \|u\|_{L^q_t,L^r_x},
$$
where $6^-<6$ (respectively $4^+>4$) is an arbitrary fixed number,
close to $6$ (respectively, to $4$). We will also write, if $I$ is
an interval and $\chi_I$ is the characteristic function of $I$,
$$
\|u\|_{S \left(I;\dot{H}^{1/2}\right)} = \left\|\chi_I u\right\|_{S
\left(\dot{H}^{1/2}\right)}.
$$

If $a$ and $b$ are two positive functions of $t$, we will write
$a=O(b)$ when there exists a constant $C>0$ (independent of $t$)
such that $a(t) \leq C b(t)$, for all $t$, and $a\approx b$ when
$a=O(b)$ and $b=O(a)$.

Throughout the paper, $C$ denotes a large positive constant and $c$
a small positive constant, that do not depend on the parameters and
may change from line to line.

\subsection*{Acknowledgements.}
T.D. was partially supported by the French ANR Grant ONDNONLIN. Part
of this work was done during S.R. visit to the University of
Cergy-Pontoise funded by the Grant ONDNONLIN. Both authors would
like to thank Frank Merle for fruitful discussions on the subject.

\section{Preliminaries}
\subsection{Properties of the ground state}
\label{SS:Q}

We recall some well-known properties of the ground state and refer
the reader to \cite{Co72}, \cite{We83}, \cite{Kw89} as well as
\cite{Ca03}, \cite[Appendix B]{Ta06BO}, \cite[\S 3]{HoRo07b} for
more details.

Consider the nonlinear elliptic equation
\begin{equation}
 \label{E:Q}
-Q+\Delta Q+|Q|^2 Q=0.
\end{equation}
The $H^1$ solutions of this equation can be enumerated by their mass
($L^2$ norm) and the minimal mass solution, $Q$, is called the
ground state. The function $Q$ is radial, smooth, positive,
exponentially decaying at infinity, and characterized as the unique
minimizer for the Gagliardo-Nirenberg inequality: if $u\in
H^1(\RR^3)$,
\begin{equation}
 \label{E:GN}
\|u\|_{4}^4 \leq C_{GN} \|\nabla u\|_{2}^3\|u\|_{2},\quad
\|Q\|_{4}^4=C_{GN} \|\nabla Q\|_{2}^3\|Q\|_{2},
\end{equation}
and
\begin{equation}
 \label{E:uniqueness}
\|u\|_{4}^4 = C_{GN}\|\nabla u\|_{2}^3\|u\|_{2}\Longrightarrow
\exists \lambda_0\in \CC,\; \exists x_0\in \RR^3\; :\;
u(x)=\lambda_0 Q(x+x_0).
\end{equation}

The above characterization of $Q$ and the concentration-compactness
principle (see \cite[Theorem I.2]{Li84IHPb}) yield:
\begin{prop}
\label{P:CC}
There exists a function $\eps(\rho)$, defined for small $\rho>0$, such that
$\lim_{\rho\rightarrow 0} \eps(\rho)=0$
and
\begin{multline*}
\forall u\in H^1, \quad \Big|\|u\|_4-\|Q\|_4\Big|+\Big|\|u\|_2 - \|Q\|_2\Big|
+ \Big|\|\nabla u\|_2-\|\nabla Q\|_2\Big|\leq \rho\\
\Longrightarrow \exists \theta_0\in \RR,\;x_0\in \RR^3,\quad \left\|
u - e^{i\theta_0}Q(\cdot-x_0)\right\|_{H^1}\leq \eps(\rho).
\end{multline*}
\end{prop}
We will make the statement of Proposition \ref{P:CC} more precise
in \S \ref{S:modulation}.

We will also need the following equalities, consequences of
Pohozhaev identities (see e.g \cite[\S 3]{HoRo07b}):
\begin{equation}
 \label{E:propQ}
\|Q\|_{4}^4 = 4 \,\|Q\|_{2}^2\quad \text{and} \quad \|\nabla
Q\|_{2}^2= 3 \|Q\|_{2}^2.
\end{equation}

\subsection{The Cauchy problem \eqref{CP}}
\label{SS:CP}
Here, we briefly recall global existence and scattering results for
\eqref{CP}, for more details see \cite{HoRo07b}. The small data
theory states that there exists a small $\ep_{sd}>0$ such that if
\begin{equation}
 \label{E:scattering_condition}
\| e^{i t \lap} u_0\|_{S([0, +\infty), \dot{H}^{1/2})} \leq
\ep_{sd},
\end{equation}
then the solution $u$ of \eqref{CP} has $T_+(u_0)=+\infty$ and
\begin{equation}
 \label{E:small_scattering}
\exists \, C>0: ~ \|u\|_{S([0, +\infty);
\dot{H}^{1/2})} \leq C \|e^{it \lap} u_0\|_{S([0, +\infty);
\dot{H}^{1/2})}.
\end{equation}
Next, if $u(t)$ is a solution which is globally defined for positive
time, then it scatters in $H^1$ as $t \to +\infty$, meaning that for
some $\phi \in H^1$, ~ $ \ds \lim_{t \to +\infty} \|u(t) - e^{i t
\lap} \phi\|_{H^1} = 0$, if it has a uniformly bounded $H^1$
norm for $t \geq 0$ and a finite Strichartz $S([0, +\infty);
\dot{H}^{1/2})$ norm. Similar statements hold for negative times.

\subsection{Gradient separation}
\label{SS:gradient_separation}
\begin{lemma}
 \label{L:gradientseparation}

Consider \eqref{CP} and suppose \eqref{E:ME-Q} holds.
\begin{enumerate}
\item \label{I:equality}
If $\|u_0\|_{2}\|\nabla u_0\|_{2} =\|Q\|_{2}\|\nabla
Q\|_{2}$, then $u=e^{it}Q$ up to the invariance of the equation.
\item
 \label{I:global}
If $\|u_0\|_{2} \|\nabla u_0\|_{2} < \|Q\|_{2}\|\nabla
Q\|_{2}$, then $u$ is globally defined and $\|u_0\|_{2} \|\nabla u (t)
\|_{2}<\|Q\|_{2} \|\nabla Q\|_{2}$ for all $t$.
\item
 \label{I:finitetime}
If $\|u_0\|_{2} \|\nabla u_0\|_{2} > \|Q\|_{2}\|\nabla
Q\|_{2}$, then $\|u_0\|_{2} \|\nabla u (t) \|_{2} >
\|Q\|_{2}\|\nabla Q\|_{2}$ for all $t$ in the domain of existence of $u$.
\end{enumerate}
\end{lemma}

\begin{proof}

Without loss of generality we can assume $M[u] = M[Q]$ and
$E[u]=E[Q]$ due to scaling: if $M[u] = \alpha M[Q]$ for some $\alpha
> 0$, then define
$\tilde u (x,t) = \alpha \, u (\alpha x, \alpha^2 t)$ which is also
a solution of \eqref{CP}, and observe that $M[\tilde u] =
\tfrac1{\alpha} M[u] = M[Q]$ and also $E[\tilde u] = \alpha E[u] =
\alpha \tfrac{M[Q]E[Q]}{M[u]} = E[Q]$.

Case \eqref{I:equality} is given by the variational characterization
\eqref{E:uniqueness} of $Q$  and the uniqueness of solutions of
\eqref{CP}.

For case \eqref{I:global} we show that if $\|\nabla u(t)\|_{2} < \|
\nabla Q\|_{2}$ holds for $t=0$, then it does for all $t$. To the
contrary, suppose (by continuity) there exists $t_1$ such that $\|
\nabla u(t_1) \|_{2} = \| \nabla Q\|_{2}$ . Then by case
\eqref{I:equality} with the initial condition at $t=t_1$, the
equality holds for all times $t$ contradicting the condition at
$t=0$. Hence, such $t_1$ does not exist and the gradient of $u(t)$
is bounded as claimed. By the finite blow-up criterion $u$ is
globally defined.

Case \eqref{I:finitetime} is similar to case \eqref{I:global}.
\end{proof}

\subsection{Spectral properties of the linearized operator}
\label{SS:linearized}

Consider a solution $u$ of \eqref{CP} close to $e^{it}Q$ and write
$u$ as
$$
u(x,t) = e^{it}(Q(x)+h(x,t)).
$$
Let $h_1=\re h$ and $h_2=\im h$. We will often identify $\CC$ and
$\RR^2$ and consider $h = h_1+i h_2 \in \CC$ as an element
$\begin{pmatrix} h_1\\h_2\end{pmatrix}$ of $\RR^2$. Note that $h$ is
a solution of the equation
\begin{equation}
 \label{lin.h}
\partial_t h +\LLL h=R(h),\quad \LLL: =\begin{pmatrix} 0& -L_-\\
L_+ &0\end{pmatrix},
\end{equation}
where the self-adjoint operators $L_+$ and $L_-$ and the remainder
$R(h)$ are defined by
\begin{gather}
\label{defL+L-}
L_+h_1:=-\Delta h_1+h_1-3 Q^2 h_1,\quad L_-h_2:=-\Delta h_2+h_2-Q^2 h_2,\\
\label{defR} R(h):= i Q\,(2|h|^2+h^2)+i|h|^2h.
\end{gather}
The spectral properties of the operator $\LLL$ are well known and
for the following proposition we refer to \cite[Theorem 3.1 and Corollary 3.1]{Gr90}
and  \cite[Proposition 2.8]{We85}.
\begin{prop}
\label{P:spectralL} Let $\sigma(\LLL)$ be the spectrum of the
operator $\LLL$, defined on $L^2(\RR^3)\times L^2(\RR^3)$
and let
$\sigma_{\ess}(\LLL)$ be its essential spectrum. Then
$$
\sigma_{\ess}(\LLL)=\left\{i\xi:\; \xi\in \RR, \; |\xi|\geq 1\right\},\quad \sigma(\LLL)\cap
\RR=\{-e_0,0,+e_0\} \quad \mbox{with} ~~ e_0>0.
$$
Furthermore, $e_0$ and $-e_0$ are simple eigenvalues of $\LLL$ with
eigenfunctions $\YYY_+$ and $\YYY_-=\overline{\YYY}_+$ in $\SSS$,
and the null-space of $\LLL$ is spanned by the four vectors
$\partial_{x_j}Q$, $j=1,2,3$ and $iQ$.
\end{prop}
\begin{remark}
\label{R:Y1Y2}
Let $\YYY_1=\re \YYY_+=\re\YYY_-$ and $\YYY_2=\im \YYY_+=-\im \YYY_-$. Then
$$L_+ \YYY_1 = e_0 \YYY_2\quad \text{and}\quad L_-\YYY_2 = - e_0 \YYY_1.$$
\end{remark}
\begin{remark}
 \label{R:null.space}
Proposition \ref{P:spectralL} implies that the null-space of $L_+$
is spanned by $\partial_{x_1}Q$, $\partial_{x_2}Q$ and
$\partial_{x_3}Q$ and the null-space of $L_-$ is spanned by $Q$.
\end{remark}
\begin{remark}
 \label{R:L-positive}
It also follows from Proposition 2.8 of \cite{We85} that $\int
(L_-f) f\geq 0$ for all real-valued $h$ in $H^1$. Together with the
preceding remark, we get
$$
\forall f\in H^1\setminus\{ \lambda Q,\;\lambda \in \RR\}, \quad \int (L_- f)f>0.
$$
\end{remark}

Consider the linearized equation $\partial_t h+\LLL h = 0$. Multiply
by $\overline{i \partial_t h}$ and take the real part to obtain
\begin{equation}
 \label{E:prePhi}
\partial_t \int (L_+ h_1)\,h_1 + \partial_t \int (L_- h_2)\,h_2 = 0.
\end{equation}
Define $\Phi$, a linearized energy, by
\begin{equation}
 \label{E:Phi}
\Phi(h):=\frac 12 \int |h|^2+\frac 12 \int |\nabla h|^2-\frac 12
\int Q^2 (3h_1^2+h_2^2)=\frac 12 \int (L_+h_1)\,h_1+\frac 12 \int
(L_-h_2)\,h_2.
\end{equation}
From \eqref{E:prePhi} it follows that $\Phi$ is conserved for
solutions of the linearized equation $\partial_t h+\LLL h=0$. By
explicit calculation (see the beginning of Appendix
\ref{A:coercivity} for the details),
\begin{equation}
 \label{BoundPhi}
E[Q+h] = E[Q] ~~\text{and}~~ M[Q+h]=M[Q]\Longrightarrow
\Phi(h)=\int Q |h|^2 h_1+\frac14 \int |h|^4,
\end{equation}
which shows that $|\Phi(h)|\leq c \, \|h\|_{4}^3$ for a threshold
solution $u=e^{it}(Q+h)$ of \eqref{CP} which is close to $e^{it}Q$.
To take advantage of this, we next study the
sign of $\Phi(h)$.

We denote by $B(g,h)$ the bilinear symmetric form associated to
$\Phi$, i.e., for $g,h\in H^1(\RR^3)$
\begin{equation}
 \label{bilinear}
B(g,h)=\frac 12 \int (L_+ g_1)\,h_1+\frac 12\int (L_- g_2)\,h_2.
\end{equation}
By Remark \ref{R:null.space},
\begin{equation}
 \label{E:orthobilinear}
\forall \, h\in H^1(\RR^3), \quad
B(\partial_{x_1}Q,h)=B(\partial_{x_2}Q,h)=B(\partial_{x_3}Q,h)=B(iQ,h)=0.
\end{equation}
Furthermore, by \eqref{E:propQ}
\begin{equation}
\label{E:PhiQ<0}
\Phi(Q) = \frac 12 \int |Q|^2+\frac 12 \int |\nabla Q|^2-\frac 32
\int Q^4=-4\int Q^2<0.
\end{equation}

Together with \eqref{E:orthobilinear} we get
$$
\forall \, h \in
\vect\{\partial_{x_1}Q,\partial_{x_2}Q,\partial_{x_3}Q,iQ,Q\}, \quad
\Phi(h)\leq 0.
$$
We next define two subspaces of $H^1$ where $\Phi$ is positive.
Consider the following orthogonality relations:
\begin{gather}
 \label{ortho1}
\int (\partial_{x_1}Q) h_1 = \int (\partial_{x_2}Q) h_1
=\int (\partial_{x_3}Q) h_1=\int Q h_2=0,\\
 \label{ortho2}
\int \Delta Q h_1=0,\\
 \tag{\ref{ortho2}'}\label{orthoef}
\int \YYY_1 \, h_2 = \int \YYY_2 \, h_1=0.
\end{gather}

Let $G_{\bot}$ be the set of $h\in H^1$ satisfying the orthogonality
relations \eqref{ortho1} and \eqref{ortho2} and $G_{\bot}'$ the set
of $h \in H^1$ satisfying \eqref{ortho1} and \eqref{orthoef}. Then
\begin{prop}[Coercivity of $\Phi$]
 \label{P:Coercivity}
There exists a constant $c>0$ such that
\begin{equation}
\label{coercivity}
\forall h\in G_{\bot}\cup G'_{\bot}, \quad \Phi(h)\geq c\|h\|_{H^1}^2.
\end{equation}
\end{prop}
Proposition \ref{P:Coercivity} is proven in appendix \ref{A:coercivity}.
Observe that as a consequence of Proposition \ref{P:Coercivity},
\begin{equation}
\label{E:nonzero}
\int \Delta Q\YYY_1\neq 0.
\end{equation}
Indeed, assume $\int \Delta Q\YYY_1=0$. Then by \eqref{E:Q} and
direct computations, $\int L_+ Q \YYY_1=0$. By Remark \ref{R:Y1Y2}
we obtain $\int Q\YYY_2=0$ which shows that $Q$ is in $G_{\bot}'$,
contradicting \eqref{E:PhiQ<0}.

\begin{remark}
In \cite{We85} Weinstein gives a sharp description of the semi-group
$e^{-t\LLL}$ for the mass-subcritical and mass-critical focusing NLS
equations. In both cases, one may decompose $H^1$ as a direct sum
$S\oplus M$, where $S$ and $M$ are stable by the flow of
$e^{-t\LLL}$, $S$ is finite dimensional and contains the eigenfunctions 
of $\LLL$, and $\Phi$ is equivalent to the $H^1$ norm on $M$, which
implies that $e^{-t\LLL}h_0$ is  bounded in $H^1$ if $h_0\in M$. It
is not clear whether such a convenient decomposition exists for the
mass-supercritical NLS equation. Note that the vector space
$G_{\bot}'$, which will play the roles of $M$ in the sequel, is not
invariant by the flow of the semi-group $e^{-t\LLL}$. However,
Proposition \ref{P:Coercivity} is sufficient for our needs, namely
the description of the dynamics of exponentially decaying solutions
of the linearized equation (see Subsection \ref{SS:expo}).
\end{remark}

\section{Existence of special solutions}
\label{S:spectralsolution}

The aim of this section is to construct the solutions $Q^+$ and
$Q^-$ of Theorem \ref{T:existence}. Namely, we will show:
\begin{prop}
\label{P:UA} Let $A\in \RR$. If $t_0=t_0(A)>0$ is large enough, then
there exists a radial solution
$$
U^A \in C^{\infty} \bigg([t_0,+\infty), H^{\infty}\bigg)
$$
of \eqref{CP} such that
\begin{equation}
 \label{E:smallUa}
\forall b\in \RR,\;\exists \, C>0 : ~ \forall t\geq t_0 \quad \text{we~have} \quad
\left\|U^A(t)-e^{it}Q-A e^{(i-e_0)t}\YYY_+\right\|_{H^b}\leq
Ce^{-2e_0t}.
\end{equation}
\end{prop}
\begin{remark}
 \label{R:U+U-}
Note that by \eqref{E:smallUa}, the conservation of mass and energy, we have
$$
M[U^A]=M[Q],\text{ and } E[U^A]=E[Q].
$$
Furthermore, again by \eqref{E:smallUa},
$$
\left\|\nabla U^A(t)\right\|^2_2 = \left\|\nabla Q \right \|^2_2
+2Ae^{-e_0t} \int \nabla Q\cdot \nabla \YYY_1+O(e^{-2e_0t}),\quad
t\rightarrow+\infty.
$$
By \eqref{E:nonzero}, replacing $\YYY_+$ by $-\YYY_+$ if necessary,
we may assume
$$
\int \nabla Q\nabla \YYY_1>0,
$$
which shows that $\left\|\nabla U^A(t)\right\|^2_2-\left\|\nabla
Q\right\|^2_2$ has the sign of $A$ for large positive time. Thus, by
Lemma \ref{L:gradientseparation}, $\left\|\nabla
U^A(t_0)\right\|^2_2-\left\|\nabla Q\right\|^2_2$ has the sign of
$A$. Letting
$$
Q^+(x,t)=e^{-it_0}U^{+1}(x,t+t_0), \quad
Q^-(x,t)=e^{-it_0}U^{-1}(x,t+t_0),
$$
we get two solutions satisfying
$$
E[Q]=E[Q^{\pm}],\quad M[Q]=M[Q^{\pm}],\quad \|\nabla Q^+(0)\|^2_2
>\|\nabla Q\|^2_2,\quad \|\nabla Q^-(0)\|^2_2<\|\nabla Q\|^2_2
$$
and
$$
\left\|Q^{\pm}-e^{it}Q\right\|_{H^1}\leq Ce^{-e_0t},\quad t\geq 0.
$$
To conclude the proof of Theorem \ref{T:existence}, it remains to
specify the behavior of $Q^+$ and $Q^-$ for negative $t$, which we
will do in Remark \ref{R:U+} and \S \ref{SS:conclu_existence}.
\end{remark}

\begin{remark}
We will see in \S \ref{SS:end_classification} that all solutions $U^A$,
$A>0$ (respectively $A<0$) are equal to $Q^+$ (respectively $Q^-$)
up to a translation in time and a multiplication by a complex number of modulus $1$.
\end{remark}

The proof of Proposition \ref{P:UA} is similar to the one of
Proposition 6.1 in \cite{DuMe07}. We start with the construction of
a family of approximate solutions to \eqref{CP} that satisfy
\eqref{E:smallUa}, and then prove the existence of $U^A$ by a fixed
point argument around an approximate solution.

\subsection{A family of approximate solutions}
\begin{prop}
 \label{P:Approximate}
Let $A\in \RR$. There exists a sequence $(\ZZZ_j^A)_{j\geq 1}$ of
functions in $\SSS$ such that $\ZZZ_1^A=A\YYY_+$ and, if
$k\geq 1$ and $\VVV_k^A:=\sum_{j=1}^k e^{-je_0t}\ZZZ_j^A$, then as
$t\rightarrow +\infty$ we have
\begin{equation}
 \label{eqVk}
\partial_t \VVV_k^A +\LLL \VVV_k^A=R(\VVV_k^A) + O\left(e^{-(k+1)e_0t}\right)
\text{ in }\SSS,
\end{equation}
where the linear operator $\LLL$ and the nonlinear term $R$ are defined in \eqref{lin.h}.
\end{prop}

\begin{remark}
Let $U_k^A:=e^{it}(Q+\VVV_k^A)$. Then $U_k^A$ is an approximate
solution of \eqref{CP} for large $t$ and satisfies
\eqref{E:smallUa}. Indeed, as $t \rightarrow +\infty$, we have
$$
i\partial_t U_k^A+\Delta U_k^A+ \left|U_k^A\right|^2 U_k^A =
O\left(e^{-(k+1)e_0t} \right) \text{ in }\SSS.
$$
\end{remark}

\begin{proof}[Proof of Proposition \ref{P:Approximate}]
We prove this proposition by induction. For brevity, we omit the
superscript $A$.

Define $\ZZZ_1:=A \YYY_+$ and $\VVV_1:=e^{-e_0t}\ZZZ_1$. Then
$$
\partial_t \VVV _1 + \LLL \VVV_1-R(\VVV_1)= -R(\VVV_1)
= -R \left(e^{-e_0t}\ZZZ_1\right),
$$
which yields \eqref{eqVk} for $k=1$.

Let $k\geq 1$ and assume that $\ZZZ_1,\ldots,\ZZZ_k$ are known with
the corresponding $\VVV_k$ satisfying \eqref{eqVk}. Expand the
expression of $R(\VVV_k)$ by using \eqref{defR}, and observe that
$R(\VVV_k)$ is of the form $\sum_{j=2}^{3k} e^{-je_0t}f_{jk}$ with
$f_{jk}$'s being in $\SSS$. Thus, by \eqref{eqVk}, there exists
$\UUU_{k+1}\in \SSS$ such that, as $t\rightarrow +\infty$, we have
$$
\partial_t \VVV_k + \LLL \VVV_k = R(\VVV_k) + e^{-(k+1)e_0 t} \UUU_{k+1}
+ O \left(e^{-(k+2)e_0t} \right) \text{ in } \SSS.
$$
By Proposition \ref{P:spectralL}, $(k+1)e_0$ is not in the spectrum
of $\LLL$. Define $\ZZZ_{k+1}:= -\left(\LLL-(k+1)e_0\right)^{-1}
\UUU_{k+1}$. It is classical that $\ZZZ_{k+1}\in \SSS$ (e.g., see
\cite[Appendix 7.2.2]{DuMe07} for an elementary proof in a similar
setting). Then we have
\begin{equation}
 \label{eqVk1}
\partial_t\left(\VVV_k +e^{-(k+1)e_0t}\ZZZ_{k+1}\right)
+ \LLL \left(\VVV_k+e^{-(k+1)e_0t}\ZZZ_{k+1}\right)\\
=R(\VVV_k)+ O\left(e^{-(k+2)e_0t} \right)\text{ in } \SSS.
\end{equation}
Denote $\VVV_{k+1}:= \VVV_k+e^{-(k+1)e_0 t} \ZZZ_{k+1}$. By
\eqref{eqVk1}, $\VVV_{k+1}$ satisfies, as $t\rightarrow +\infty$,
\begin{equation}
 \label{eqVk+1}
\partial_t \VVV_{k+1}+\LLL \VVV_{k+1} - R(\VVV_{k+1})=
R(\VVV_{k})-R(\VVV_{k+1})+O\left(e^{-(k+2)e_0 t}\right)\text{ in }\SSS.
\end{equation}
Since we have, as $t\rightarrow +\infty$,
$\VVV_j=O\left(e^{-e_0t}\right)$ in $\SSS$ for $j=k,k+1$ and
$\VVV_k-\VVV_{k+1}=O(e^{-(k+1)e_0t})$ in $\SSS$, we obtain (using
the explicit expression of $R$),
\begin{equation*}
R(\VVV_{k})-R(\VVV_{k+1}) = O\left(e^{-(k+2)e_0 t}\right) \text{ in
}\SSS,
\end{equation*}
as $t\rightarrow +\infty$ which shows, in view of \eqref{eqVk+1}, the desired estimate
\eqref{eqVk} for $k+1$. This completes the proof.
\end{proof}

\subsection{Construction of special solutions}

Next we prove Proposition \ref{P:UA}. We will construct a solution
$U^A$ such that there exists $t_0 \in \cR$
\begin{equation}
 \label{E:asymptoticUa}
\forall \, b \in \cR,\; \exists \, C>0:  \; \forall \, t\geq t_0
~\text{and} ~~ \forall \, k \in \cN, ~~ \left\| U^A(t) - e^{it}(Q +
\VVV_k^A(t) ) \right\|_{H^b} \leq C e^{-(k+1)e_0t}.
\end{equation}
Let $b>3/2$ and write
$$
U^A=e^{it}\left(Q+h^A\right).
$$
First, by a fixed point argument we construct a
solution of \eqref{lin.h} $h^A \in C^{0}([t_k,+\infty),H^b)$ for $k$
and $t_k$ large and such that
\begin{equation}
 \label{estimh}
\forall \, t\geq t_k, ~~ \|h^A(t)-\VVV_k^A(t)\|_{H^b} \leq
e^{-(k+\frac12)e_0t}.
\end{equation}
Next, we show  by uniqueness arguments that $h^A$ does not depend on
$b$ and $k$. Estimate \eqref{E:asymptoticUa} will follow from
\eqref{estimh}. For brevity we again omit the superscript $A$.

\medskip

\noindent\emph{Step 1. Reduction to a fixed point problem.}
The equation \eqref{lin.h} may be written as a Schr\"odinger equation
\begin{equation}
 \label{eqhNLS}
i \partial_t h+\Delta h-h=-S(h),\quad
S(h):=2Q^2h+Q^2\overline{h}+2Q|h|^2+Qh^2+|h|^2h.
\end{equation}
For $k \in \cN$ define
\begin{equation}
 \label{eqepsk}
\eps_k(t)=i\partial_t \VVV_k+\Delta\VVV_k- \VVV_k+S(\VVV_k).
\end{equation}
By Proposition \ref{P:Approximate}, as $t \to +\infty$,
\begin{equation}
 \label{eqepsasymptotic}
\eps_k(t) = O(e^{-(k+1)e_0t})\text{ in }\SSS.
\end{equation}
Let $v:=h-\VVV_k$ and subtract \eqref{eqepsk} from \eqref{eqhNLS} to
obtain
\begin{equation}
 \label{eqvfp}
i\partial_t v+\Delta v-v=-S(\VVV_k+v)+S(\VVV_k)-\eps_k.
\end{equation}
The solution of \eqref{eqvfp} is given by the equation
\begin{gather}
\label{fxpt}
v(t)=\MMM(v)(t),\\
\notag
\text{where }\MMM(v)(t):=-i\int_t^{+\infty} e^{i(t-s)
(\Delta-1)}\Big[S(\VVV_k(s)+v(s))- S(\VVV_k(s))+ \eps_k(s)\Big] \,ds.
\end{gather}
Note that \eqref{estimh} is equivalent to $\ds \|v(t)\|_{H^b} \leq
e^{-\left(k + \frac 12\right)e_0t}$ for $t\geq t_k.$ Thus, we must
show that $\MMM$ is a contraction on $B$ defined by
\begin{align*}
B&=B(t_k,k,b):=\Big\{v\in E,\; \|v\|_{E}\leq 1\Big\},\\
E&=E(t_k,k,b):=\Big\{v\in C^{0}([t_k,+\infty),H^b),\; \|v\|_{E}=\sup_{t\geq t_k}
e^{\left(k+\frac 12\right)e_0t}\|v(t)\|_{H^b}<\infty\Big\}.
\end{align*}
This is the object of the following step.

\medskip

\noindent\emph{Step 2. Contraction argument.} We show that $\MMM$ is
a contraction on $B$ for $b>3/2$, and $k$ and $t_k$ sufficiently
large\footnote{Note that the condition $b>3/2$ is not restrictive:
if \eqref{E:asymptoticUa} is shown for some $b_0$, it follows for
all $b\leq b_0$.}

Throughout this proof, we denote by $C$ a constant depending only on
$b$, and $C_k$ a constant depending on $b$ and $k$. Both constants
may change from line to line. Note that $H^b$ is closed under
multiplication and conjugation for $b>\frac 32$. In view of the
identities
\begin{gather*}
F^2-G^2=(F-G)(F+G),\quad |F|^2-|G|^2=\re\big((\overline{F}
-\overline{G})(F+G)\big),\\
|F|^2F-|G|^2G=F\re\big((\overline{F}-\overline{G})(F+G)\big)+|G|^2(F-G),
\end{gather*}
we obtain that for $F,G\in H^b$ there exists a constant $C_0>0$ such
that
\begin{equation}
 \label{ineqS}
\big\|S(F)-S(G)\big\|_{H^b} \leq C_0\big\|F-G\big\|_{H^b}
\left(1+\big\|F\big\|_{H^b}^2+\big\|G\|_{H^b}^2\right).
\end{equation}
Let $v\in B$. Observe that for all $t\in \RR$, $e^{it(\Delta-1)}$ is
an isometry of $H^b$. By the definition of $\MMM$, and applying the
bound \eqref{eqepsasymptotic} on $\eps_k$ and the estimate
\eqref{ineqS}, we get
\begin{multline}
 \label{ineqM1}
\forall \,t\geq t_k,\quad \|\MMM(v(t))\|_{H^b}\leq C\int_t^{+\infty}
\|v(s)\|_{H^b}\left(1+\big\|\VVV_k(s)\big\|_{H^b}^2
+ \big\|v(s)\|_{H^b}^2\right)\,ds\\
+ C_k\int_t^{+\infty} e^{-(k+1)e_0 s} \,ds.
\end{multline}
By the construction of $\VVV_k$, $\|\VVV_k(s)\|_{H^b}\leq
C_ke^{-e_0s}$. Furthermore, since $v\in B$, $\|v(s)\|_{H^b}\leq
e^{-(k+\frac 12)e_0s}$. Hence,
\begin{multline*}
\forall \, t\geq t_k,\quad \int_t^{+\infty} \|v(s)\|_{H^b} \left(1
+\big\|\VVV_k(s)\big\|_{H^b}^2+\big\|v(s)\|_{H^b}^2\right)\,ds\\
\leq 2\int_t^{+\infty} \left(e^{-\left(k+\frac 12\right)e_0s} +
C_ke^{-\left(k+\frac 52\right)e_0s}\right) ds \leq
\left(\frac{2}{\left(k + \frac 12\right)e_0} + C_k e^{-2e_0t}\right)
e^{-\left(k+\frac 12\right)e_0t}.
\end{multline*}
Therefore, $\MMM(v)\in E$, and from \eqref{ineqM1} we obtain
$$
\|\MMM(v)\|_E\leq \frac{C}{k+\frac 12}+C_ke^{-\frac 12 e_0 t_k}.
$$
First, choose $k$ so that $\frac{C}{k+\frac 12}\leq \frac 12$; next,
take $t_k$ such that $C_ke^{-\frac 12 e_0 t_k}\leq \frac 12$. Then
$\MMM$ maps $B=B(t_k,k,b)$ into itself.

It remains to show that $\MMM$ is a contraction. Let $v,w\in B$. By
the definition of $\MMM$ and \eqref{ineqS}, we have
\begin{align*}
\big\| \MMM(v)(t)-\MMM(w)(t)\big\|_{H^b}&\leq \int_t^{+\infty}
\big\|S(\VVV_k(s)+v(s))-S(\VVV_k(s)+w(t))\big\|_{H^b}ds\\
&\leq \int_t^{+\infty} \left(C+C_ke^{-2e_0s}\right)\|v(s)-w(s)\|_{H^b}ds\\
&\leq \left(\frac{C}{k+\frac 12}+C_k e^{-2e_0t_k}\right)
e^{-\left(k+\frac 12\right)t} \|v-w\|_E.
\end{align*}
Choosing if necessary a larger $k$, then a larger $t_k$, we may
assume that $\frac{C}{k+\frac 12}<\frac 12$ and $C_ke^{-2e_0t_k}\leq
\frac 12$, showing that $\MMM$ is a contraction on $B$. Hence, Step
2 is complete.
\medskip

\noindent\emph{Step 3. End of the proof.}

By the preceding step with $b=2$, there exists $k_0$ and $t_0$ such
that there exists a unique solution $U^A$ of \eqref{CP} satisfying
$U^A\in C^0([t_0,+\infty), H^2)$ and for all $t\geq t_0$
\begin{equation}
 \label{propUA}
\left\|U^A(t)-e^{it}Q-e^{it} \VVV_{k_0}^A(t)\right\|_{H^2} \leq
e^{-(k_0+\frac{1}{2}) e_0 t}.
\end{equation}
Note that the fixed point argument still holds taking a larger
$t_0$, so that the uniqueness remains valid, for any $t_0'\geq t_0$,
in the class of solutions of \eqref{CP} in $C^0([t_0',+\infty),
H^2)$ satisfying \eqref{propUA} for all $t\geq t_0'$.

Next, we show that $U^A\in C^{\infty}([t_0,+\infty),H^b)$ and that
\eqref{E:asymptoticUa} holds for any $b \in \cR$.

Since $U^A$ is a solution of \eqref{CP}, it is sufficient to show
that $U^A \in C^{0}([t_0,+\infty),H^b)$ for any $b$; the smoothness
in time will follow from \eqref{CP} and Sobolev embeddings. Let
$b\geq 2$. By Step 2, if $k_1$ is large enough, there exists $t_1$
and $\tU^A\in C^{0}([t_1,+\infty),H^b)$ such that
$$
\forall \, t\geq t_1, ~~
\left\|\tU^A-e^{it}Q-e^{it}\VVV_{k_1}^A(t)\right\|_{H^b} \leq
e^{-(k_1+\frac{1}{2}) e_0 t}.
$$
Of course, we may choose $k_1\geq k_0+1$. Using that by the
construction of $\VVV_k^A$,
$\left\|\VVV_{k_1}^A-\VVV_{k_0}^A\right\|_{H^b}\leq C \,
e^{-(k_0+1)e_0 t}$, we get
\begin{equation}
 \label{E:asymptotictU}
\forall \, t\geq t_1, ~~
\left\|\tU^A-e^{it}Q-e^{it}\VVV_{k_0}^A(t)\right\|_{H^b}\leq
e^{-(k_1+\frac{1}{2}) e_0 t}+C e^{-(k_0+1)e_0t} \leq
Ce^{-(k_0+1)e_0t}.
\end{equation}
In particular, $\tU^A$ satisfies \eqref{propUA} for large $t$. By
the uniqueness in the fixed point argument, $U^A=\tU^A$, which shows
that $U^A\in C^0([t_1,+\infty),H^b)$, $b \geq 2$. By the persistence
of regularity of equation \eqref{CP}, $U^A\in C^0([t_0, +\infty),
H^b)$, $b \geq 2$ (and thus for any $b \in \cR$).
Finally, we note that \eqref{E:asymptotictU} 
implies \eqref{E:asymptoticUa} with $k_0$ replaced by $k$, which
concludes the proof of Proposition \ref{P:UA}. \qed

\begin{remark}
\label{R:asymptotic} The estimate \eqref{E:asymptoticUa} yields an
asymptotic expansion of $U^A$ in terms of $e^{-e_0t}$.
\end{remark}

\section{Modulation}
\label{S:modulation}

For $u\in H^1$ define
$$
\delta(u)=\left| \int |\nabla Q|^2-\int |\nabla u|^2\right|.
$$
By Proposition \ref{P:CC}, if
\begin{equation}
 \label{E:equalityEM}
M[u]=M[Q],\quad E[u]=E[Q],
\end{equation}
and $\delta(u)$ is small enough, then there exists $\tilde{\theta}$
and $\tilde{X}$ such that $ e^{-i\tilde{\theta}}
u(\cdot+\tilde{X})=Q+\tilde{u}$ with $\|\tilde{u}\|_{H^1}\leq
\tilde{\eps}\big(\delta(u)\big)$, where
$\tilde{\eps}(\delta)\rightarrow 0$ as $\delta\rightarrow 0$. Recall
that any solution such that \eqref{E:ME-Q} holds may be rescaled to
a solution satisfying \eqref{E:equalityEM}. The goal of this section
is to choose parameters $\tilde{\theta}$ and $\tilde{X}$, when $u$
is a solution of \eqref{CP}, in order to obtain linear estimates of
these parameters and their derivatives in terms of $\delta(u)$. We
only sketch the proofs, which are very close to the ones of \cite[\S
3.2]{DuMe07}.

\begin{lemma}
 \label{lem.ortho}
There exist $\delta_0>0$ and a positive function $\eps(\delta)$
defined for $0<\delta\leq \delta_0$, which tends to $0$ when
$\delta$ tends to $0$, such that for all $u$ in $H^1(\RR^3)$
satisfying \eqref{E:equalityEM} and $\delta(u)<\delta_0$, there
exists a couple $(\sigma,X)$ in $\RR\times \RR^3$ such that
$v=e^{-i\sigma}u(\cdot+X)$ satisfies
\begin{gather}
\label{closetoQ}
\|v-Q\|_{H^1}\leq \eps(\delta),\\
\label{ortho} \im \int Q \, v=0,\quad \re \int \partial_{x_k} Q \,
v=0,\; k=1,2,3.
\end{gather}
The parameters $\sigma$ and $X$ are unique in $\RR /_{\ds 2\pi\ZZ}
\times \RR^3$, and the mapping $u\mapsto (\sigma,X)$ is $C^1$.
\end{lemma}

\begin{proof}
Let us first show the lemma when $u$ is close to $Q$ in $H^1$.
Consider the following functionals on $\RR\times \RR^3\times
H^1$:
$$
J_0: (\sigma,X,u)\mapsto \im \int e^{-i\sigma} u(x+X) \, Q , \quad
J_k: (\sigma,X,u)\mapsto \re \int e^{-i\sigma} u(x+X) \,
\partial_{x_k}Q, \; k=1,2,3.
$$
Denote $J=(J_0,J_1,J_2,J_3)$. The orthogonality conditions
\eqref{ortho} are equivalent to the condition $J(\sigma,X,u)=0$.

Note that $J(0,0,Q)=0$. Furthermore, by direct computation and using
that since $Q$ is radial,
$\int \partial_{x_k}Q\partial_{x_j}Q=0$ if $j\neq k$, one can
check that $\left(\frac{\partial J}{\partial \sigma},\frac{\partial
J}{\partial X_1},\frac{\partial J}{\partial X_2},\frac{\partial
J}{\partial X_3}\right)$ is invertible at $(0,0,Q)$. By the Implicit
Function Theorem there exists $\eps_0,\eta_0>0$ such that for $u\in
H^1$
$$
\|u-Q\|_{H^1}<\eps_0\Longrightarrow \exists! (\sigma,X),\quad
|\sigma|+|X|\leq \eta_0\text{ and }J(\sigma,X,u)=0.
$$
If $u$ is as in the Lemma, we reduce the proof to the previous case
by choosing $\tilde{\theta}$ and $\tilde{X}$ as in the introduction
of this section, so that $e^{-i\tilde{\theta}}u(\cdot+\tilde{X})$ is
close to $Q$ in $H^1$. The assertions on the uniqueness of
$(\sigma,X)$ and the regularity of the mapping $u\mapsto(\sigma,X)$
follows from the Implicit Functions Theorem.
\end{proof}

Let $u$ be a solution of \eqref{CP} satisfying \eqref{E:equalityEM}.
In the sequel we will write
$$
\delta(t):=\delta(u(t)).
$$
Let $D_{\delta_0}$ be the open set of all times $t$ in the domain of
existence of $u$ such that $\delta(t)<\delta_0$. On $D_{\delta_0}$,
by Lemma \ref{lem.ortho} we can define parameters $\sigma(t)$,
$X(t)$, which are $C^1$ functions of $t$. In the forthcoming
sections, we show under the additional hypothesis that $u$ is close
to $e^{it}Q$, up to constant modulation parameters, and thus, we
rather work with the parameters $X(t)$ and $\theta(t)=\sigma(t)- t$.
Write
\begin{gather}
\label{E:decompo-u}
e^{-i\theta(t)-it}u\left(t,x+X(t)\right)=(1+\alpha(t))Q(x)+h(t,x),\\
\notag \alpha(t)=\re \frac{e^{-i\theta-it}\int \nabla
u(t,x+X(t))\cdot \nabla Q(x)dx}{\int |\nabla Q|^2}-1.
\end{gather}
Observe that $\alpha$ is chosen such that $h$ satisfies the
orthogonality condition $\eqref{ortho2}$. By Lemma \ref{lem.ortho},
$h$ also satisfies the orthogonality conditions \eqref{ortho1}.

We next obtain a first estimate on the parameters.
\begin{lemma}
 \label{LemEstimModul1}
Let $u$ be a solution of \eqref{CP} satisfying \eqref{E:equalityEM}.
Then, taking a smaller $\delta_0$ if necessary, the following
estimates hold for $t \in D_{\delta_0}$:
\begin{equation}
 \label{EstimModul1}
|\alpha(t)|\approx \left|\int Q h_1(t)\right| \approx
\|h(t)\|_{H^1}\approx \delta(t).
\end{equation}
\end{lemma}
\begin{proof}
Let $\tilde{\delta}(t):=|\alpha(t)|+\delta(t)+\|h(t)\|_{H^1}$, which
is small when $\delta(t)$ is small (see \eqref{closetoQ}). The
equality $M[Q+\alpha Q+h]=M[u]=M[Q]$ implies $\int |\alpha
Q+h|^2+2\alpha \int Q^2+2\int Qh_1 =0$, and hence,
\begin{equation}
\label{equa1} |\alpha| =\tfrac{1}{M[Q]}\left|\int Q
h_1\right|+O\left(\tilde{\delta}^2 \right).
\end{equation}

Furthermore, by definition of $\delta$ and the orthogonality
condition \eqref{ortho2} on $h$, we obtain
$$
\delta(t)=\left| \int |\nabla(Q+\alpha Q+h)|^2-\int |\nabla
Q|^2\right| =\left| (2\alpha+\alpha^2) \int |\nabla Q|^2+\int
|\nabla h|^2\right|,
$$
which yields
\begin{equation}
\label{equa2} |\alpha|=\tfrac1{2 \|\nabla Q\|^2}\, \delta +
O\left(\tilde{\delta}^2\right).
\end{equation}
Note that the orthogonality condition $\int \nabla Q\cdot \nabla
h_1=0$ implies with equation \eqref{E:Q} that $\int Q^3 h_1=\int Q
h_1$. Thus, $B(Q,h)=-\int Qh_1$, where $B$ is as in
\eqref{bilinear}.
By \eqref{BoundPhi},
\begin{gather}
 \notag
\left|\alpha^2\Phi(Q)+\Phi(h)-2\alpha \int Qh_1\right| =|\Phi(\alpha
Q+h)| =O\left(\alpha^3+\|h\|^3_{H^1}\right),\\
 \label{Phi-estimate}
\Phi(h) = \alpha^2|\Phi(Q)| +2\alpha \, \int Q\, h_1 +
O\left(\alpha^3+\|h\|_{H^1}^3\right).
\end{gather}
By Proposition \ref{P:Coercivity} and \eqref{BoundPhi},
$\Phi(h)\approx \|h\|_{H^1}^2$. Combining this and \eqref{Phi-estimate}, we get
\begin{equation}
 \label{equa3}
\|h\|_{H^1}=O\left(|\alpha|+\left|\int Q h_1\right|+\tilde{\delta}^{3/2}\right).
\end{equation}
Substituting \eqref{equa1} into \eqref{equa3}, we get
$\|h\|_{H^1}=O\left(|\alpha|+\tilde{\delta}^{3/2}\right)$, and thus,
with \eqref{equa2}, $\tilde{\delta}=O(|\alpha|)$, which shows that
\eqref{equa1}, \eqref{equa2} and \eqref{equa3} imply
\eqref{EstimModul1}.
\end{proof}

In the sequel we denote by $'$ the derivative with respect to $t$.
\begin{lemma}[Bounds on the time-derivatives]
\label{LemEstimModul2} Under the assumption of Lemma
\ref{LemEstimModul1}, taking a smaller $\delta_0$ if necessary, we
have for $t \in D_{\delta_0}$
\begin{equation}
 |\alpha'|+|X'|+|\theta'|=O(\delta).
\end{equation}
\end{lemma}
\begin{proof}
Let $\delta^*(t)=\delta(t)+|\alpha'(t)|+|X'(t)|+|\theta'(t)|$.
The equation \eqref{CP} and Lemma \ref{LemEstimModul1} yield for
$t\in D_{\delta_0}$
\begin{equation}
 \label{equavh}
i\partial_t h+\Delta h +i\alpha' Q-\theta'Q-iX'\cdot \nabla Q
=O(\delta+\delta\delta^*) \text{ in } L^2.
\end{equation}
Note that by orthogonality relations \eqref{ortho1} and
\eqref{ortho2} on $h$, we have
\begin{equation}
 \label{orthoh'}
\forall \, t\in D_{\delta_0},\quad \re \int
\partial_{x_j} Q \, \partial_t h = 0,\; j=1,2,3,
\quad \im \int  Q \, \partial_t
h = \re \int \Delta Q \, \partial_t h = 0.
\end{equation}
Multiplying \eqref{equavh} by $Q$, integrating the real part on
$\RR^3$, we get by integration by parts (and using that by
\eqref{ortho2} $\re \int h \Delta Q=0$),
\begin{equation}
\label{theta'}
 |\theta'|=O(\delta+\delta\delta^*).
\end{equation}
Similarly, multiplying \eqref{equavh} by $\partial_{x_j} Q$,
$j\in\{1,2,3\}$,  and integrating the imaginary part on $\RR^3$, we
obtain (using that by Lemma \ref{LemEstimModul1},
$\int \Delta h\partial_{x_j}Q=O(\delta)$),
\begin{equation}
 \label{Xj'}
\left|X_j'\right|=O(\delta+\delta\delta^*),
\end{equation}
where $X=(X_1,X_2,X_3)$. Multiplying \eqref{equavh} by $\Delta Q$
and taking the imaginary part we get (noting that $\int
\partial_{x_j}Q \Delta Q=0$ and that by Lemma \ref{LemEstimModul1},
$\int \Delta h\partial_{x_j}Q=O(\delta)$),
\begin{equation}
\label{alpha'}
|\alpha'|=O(\delta+\delta\delta^*).
\end{equation}
Summing up \eqref{theta'}, \eqref{Xj'} and \eqref{alpha'}, we obtain
\begin{equation*}
\delta^*=O(\delta+\delta\delta^*),
\end{equation*}
which yields the result if $\delta_0$ is chosen small enough.
\end{proof}

We conclude this section by showing the following Lemma, needed in
the next two sections.
\begin{lemma}
 \label{L:expo}
Let $u$ be a solution of \eqref{CP} satisfying \eqref{E:equalityEM}.
Assume that $u$ is defined on $[0,+\infty)$ and that there exist
$c,C>0$ such that
\begin{equation}
\label{E:CVexp}
\forall t\geq 0,\quad \int_{t}^{+\infty}\delta(s)ds\leq Ce^{-ct}.
\end{equation}
Then there exist $\theta_0\in\RR$, $x_0\in \RR^3$, $c,C>0$ such that
$$
\|u-e^{it+i\theta_0}Q(\cdot-x_0)\|_{H^1}\leq Ce^{-ct}.
$$
\end{lemma}

\begin{proof}
First observe that \eqref{E:CVexp} implies that there exists
$\{t_n\}_{n \in \cN}$ with $t_n \to +\infty$ such that
\begin{equation}
\label{E:delta0}
\lim_{n\rightarrow+\infty} \delta(t_n)=0.
\end{equation}
 If not, there exists $\eps > 0$ such that $\delta(t) \geq
\eps$ for a.e. $t>0$, which would contradict \eqref{E:CVexp}.

Fix such $\{t_n\}_{n \in \cN}$. Next we show that $\delta(t)$ tends
to $0$ as $t$ tends to $+\infty$. If not, there exists a sequence
$\{t'_n\}_{n \in \cN}$ which tends to $+\infty$ such that
$\delta(t'_n) \geq \eps_1 > 0$ for some $\eps_1 > 0$. Now we can
choose $\{t'_n\}$, extracting subsequences from $\{t_n\}$ and/or
$\{t'_n\}$ if necessary, with the following properties:
\begin{equation}
 \label{E:absurdd0a}
\forall n,\quad t_n < t'_n
\end{equation}
\begin{equation}
 \label{E:absurdd0b}
\delta(t'_n) = \eps_1,
\end{equation}
and
\begin{equation}
 \label{E:absurdd0c}
\forall \, t \in [t_n, t'_n) \quad \delta(t)<\eps_1.
\end{equation}
On $[t_n,t'_n]$ the parameter $\alpha(t)$ is well defined. By Lemma \ref{LemEstimModul2},
$|\alpha'(t)|=O(\delta(t))$, so that by \eqref{E:CVexp},
$\int_{t_n}^{t_n'} |\alpha'(t)|dt \leq Ce^{-ct_n}$. Hence,
\begin{equation}
\label{E:alphatnt'n0} \lim_{n\rightarrow +\infty}
|\alpha(t_n)-\alpha(t_n')|=0.
\end{equation}
By Lemma \ref{LemEstimModul1}, we have $|\alpha(t)|\approx
\delta(t)$. As a consequence,  \eqref{E:delta0}
implies that $|\alpha(t_n)|$ tends to $0$, however,
\eqref{E:absurdd0b} implies that $|\alpha(t_n')|$ is bounded from
below as $t$ tends to $+\infty$. This contradicts
\eqref{E:alphatnt'n0} and shows as announced
\begin{equation}
 \label{E:d0}
\lim_{t\rightarrow +\infty} \delta(t)=0.
\end{equation}

To conclude the proof of Lemma \ref{L:expo}, in view of the
decomposition \eqref{E:decompo-u} of $u$, it is sufficient to show
that there exists $\theta_{\infty}$, $X_{\infty}$ such that
\begin{equation}
 \label{E:estimexp}
\forall t \geq 0,\quad \delta(t) + |\alpha(t)| + \|h(t)\|_{H^1} +
|X(t)-X_{\infty}| + |\theta(t)-\theta_{\infty}|\leq Ce^{-ct}.
\end{equation}
Note that by \eqref{E:d0} and the estimate $|\alpha(t)|\approx
\delta(t)$, $\alpha(t)$ tends to $0$ as $t$ tends to $+\infty$, and
thus, $\alpha(t)=-\int_t^{\infty}\alpha'(s)ds$. By \eqref{E:CVexp}
and the estimate $|\alpha'(t)|=O(\delta(t))$, we get the bound on
$\alpha(t)$ in \eqref{E:estimexp}. Since by Lemma
\ref{LemEstimModul1}, $|\alpha(t)|\approx \|h(t)\|_{H^1}\approx
\delta(t)$, we deduce the bound on $\delta$ and $h$. From Lemma
\ref{LemEstimModul2}, we get $|X'(t)|+|\theta'(t)|\leq Ce^{-ct}$.
Thus, there exist $X_{\infty}$, $\theta_{\infty}$ such that
$|X(t)-X_{\infty}|+|\theta(t)-\theta_{\infty}|\leq Ce^{-ct}$,
concluding the proof of \eqref{E:estimexp}.
\end{proof}

\section{Convergence to $Q$ in the case
$\|\nabla u_0\|_2\|u_0\|_2>\|\nabla Q\|_2\|Q\|_2$}
\label{S:supercritical}

In this section, we show the following proposition, which is the
first step in the proof of case \eqref{case_super} of Theorem
\ref{T:classification}.
\begin{prop}
\label{P:supercrit}
Consider a solution $u$ of \eqref{CP} such that
\begin{gather}
 \label{E:equality}
M[u]=M[Q],\quad E[u]=E[Q],\\
\label{E:supercrit}
\|\nabla u_0\|_2>\|\nabla Q\|_2,
\end{gather}
which is globally defined for positive times. Assume furthermore that
$u_0$ is either of finite variance, i.e.,
\begin{equation}
 \label{E:finitevar}
\int |x|^2 |u_0(x)|^2 \, dx < \infty,
\end{equation}
or radial.
Then there exists $\theta_0\in\RR$, $x_0\in\RR^3$, $c,C>0$ such that
$$
\|u-e^{it+i\theta_0}Q(\cdot-x_0)\|_{H^1}\leq Ce^{-ct}.
$$
Moreover, the negative time of existence of $u$ is finite.
\end{prop}

\begin{remark}
 \label{R:U+}
The last statement of Proposition \ref{P:supercrit} shows that the
radial solution $Q^+$ constructed in Remark \ref{R:U+U-} has finite
negative time of existence.
\end{remark}

\subsection{Finite variance solutions}
\label{SS:supercrit-finitevariance}

In the finite variance case, Proposition \ref{P:supercrit} relies on
the following lemma.
\begin{lemma}
 \label{L:supercrit-finitevariance}
Let $u$ be a solution of \eqref{CP} satisfying \eqref{E:equality},
\eqref{E:supercrit}, \eqref{E:finitevar} and
\begin{gather}
 \label{as.supercrit2}
T_+(u_0)=+\infty.
\end{gather}
Then for all $t$ in the interval of existence of $u$,
\begin{equation}
 \label{super.positive}
\im \int  x \cdot \nabla u(x,t) \, \ubar(x,t) \, dx>0,
\end{equation}
and there exist $c,C>0$ such that
\begin{equation}
\label{super.expdecay} \forall ~ t\geq 0,\quad \int_t^{\infty}
\delta(s)\, ds \leq C\,e^{-c\,t}.
\end{equation}
\end{lemma}
Let us first assume Lemma \ref{L:supercrit-finitevariance} and prove
Proposition \ref{P:supercrit} in the finite variance case.

By \eqref{super.positive}, $\im \int  x \cdot \nabla u(x,t) \,
\ubar(x,t) \, dx >0$ for all $t$ in the interval of existence of
$u$. Now assume that $u$ is also globally defined for negative
times, and consider the function $v(x,t)=\overline{u}(x,-t)$. Then
$v$ is a solution of \eqref{CP} satisfying the assumptions of Lemma
\ref{L:supercrit-finitevariance}. Thus, by \eqref{super.positive},
for all $t$ in the domain of existence of $u$,
$$
0< \im \int  x \cdot \nabla v(x,-t) \, \vbar(x,-t)=-\im \int  x
\cdot \nabla u(x,t) \, \ubar(x,t) \, dx,
$$
which contradicts \eqref{super.positive}. Hence, the negative time
of existence of $u$ is finite. The other assertion of Proposition
\ref{P:supercrit} follows from \eqref{super.expdecay} and Lemma \ref{L:expo}.

To finish Subsection \ref{SS:supercrit-finitevariance}, we must
prove Lemma \ref{L:supercrit-finitevariance}.

\begin{proof}[Proof of Lemma \ref{L:supercrit-finitevariance}]
Let $y(t):=\int |x|^2 |u(x,t)|^2 \, dx$. Then by equation \eqref{CP} and direct
computations, we have, for all $t$ in the interval of existence of
the solution $u$ 
\begin{equation}
 \label{y'bis}
y'(t) = 4 \im \int x \cdot \nabla u \, \ubar \, dx.
\end{equation}
Furthermore, $y''(t)=8\int |\nabla u|^2-6\int |u|^4$. By
\eqref{E:propQ}, $\|Q\|_4^4=\frac 43 \|\nabla Q\|_2^2$. Using that
$E[u]=E[Q]$, we get $8\int |\nabla u|^2-6\int |u|^4=4\left(\|\nabla
Q\|^2_2-\|\nabla u\|^2_2\right)$, and thus,
\begin{equation}
 \label{y''bis}
y''(t) =-4\,\delta(t) <0.
\end{equation}

\noindent\emph{Step 1.} We first show \eqref{super.positive}, which is equivalent to
\begin{equation}
 \label{y'>0}
\quad y'(t)>0.
\end{equation}
If not, there exists $t_1$ 
such that $y'(t_1)\leq 0$. By
\eqref{y''bis}, if $t_0 > t_1$ we obtain
\begin{equation}
\label{absurd.y'<0}
\forall t\geq t_0,\quad y'(t)\leq y'(t_0)<0.
\end{equation}
As $T_+(u_0)=+\infty$, this shows that $y(t)<0$ for large $t$,
yielding a contradiction.
\medskip

\noindent\emph{Step 2. End of the proof}

We first note that
\begin{equation}
\label{super.bound}
\left(y'(t)\right)^2\leq C \,y(t)\left(y''(t)\right)^2.
\end{equation}

Indeed this is an immediate consequence of \eqref{y'bis},
\eqref{y''bis} and the following claim, in the spirit of
\cite[Lemma 2.1]{Ba04}, proven in Appendix \ref{A:claimval}.
\begin{claim}
 \label{C:val}
Let $\varphi \in C^{1}(\RR^3)$ and $f\in H^1(\RR^3)$. Assume that
$\int |f|^2 |\nabla \varphi|^2$ is finite and
\begin{equation}
\label{hyp.val} \|f\|_2=\|Q\|_2,\quad E[f]=E[Q].
\end{equation}
Then
\begin{equation}
\label{CauSchVal} \left\vert \im \int (\nabla \varphi \cdot \nabla
f) \, \fbar \, \right\vert^2 \leq C \, \delta^2(f) \int |\nabla
\varphi|^2 |f|^2 .
\end{equation}
\end{claim}
Taking $\varphi(x)=|x|^2$ in the above Claim,  we obtain
\eqref{super.bound}.
\smallskip

Now, for all $t$ in the interval of existence of $u$ we get
$y'(t)>0$ and $y''(t)<0$ (see \eqref{y'>0} and \eqref{y''bis}).
Thus,
\begin{equation}
\label{E:diffineq}
\quad \frac{y'(t)}{\sqrt{y(t)}}\leq -C y''(t).
\end{equation}
Integrating between $0$ and $t\geq 0$, we get
$$
\sqrt{y(t)}-\sqrt{y(0)} \leq
-C(y'(t)-y'(0)) \leq Cy'(0).
$$
This shows that $y(t)$ is bounded for $t \geq 0$. From
\eqref{E:diffineq} we deduce
$$
\quad y'(t)\leq -Cy''(t),
$$
which shows
$$
y'(t)\leq Ce^{-ct}.
$$

Now
$$
y'(t)=-\int_t^{\infty} y''(s) ds=4\int_t^{+\infty} \delta(s) ds.
$$
Hence, we obtain \eqref{super.expdecay}, concluding the proof of
Lemma \ref{L:supercrit-finitevariance}.
\end{proof}

\subsection{Radial solutions.}
Assume that $u$ is radial, satisfies \eqref{E:equality} and
\eqref{E:supercrit}, and that it is globally defined for positive
time. We will show that $u$ has finite variance, which will yield
Proposition \ref{P:supercrit} in the case of radial solutions also.

Let $\varphi$ be a radial function such that
$$
0\leq \varphi(r),\quad \varphi''(r)\leq 2,\quad 0\leq r\leq 1
\Longrightarrow \varphi(r)=r^2,\quad 2\leq r\Longrightarrow \varphi(r)=0.
$$
Consider the localized variance
$$
y_R(t)=\int R^2\varphi\left(\frac xR\right)|u(x,t)|^2dx.
$$
We know that \eqref{E:equality} implies $8\int |\nabla u|^2-6\int
|u|^4=4\left(\|\nabla Q\|^2_2-\|\nabla u\|^2_2\right)$ (see the
computation before \eqref{y''bis}). By explicit calculations,
\eqref{CP}, \eqref{E:equality} and the radiality of $u$ and
$\varphi$ we get
\begin{gather}
 \label{E:y'Rsupercrit}
y'_R(t)=2R\im \int\overline{u}\,\nabla\varphi\left(\frac xR\right) \cdot\nabla u,\\
 \label{E:y''Rsupercrit}
y''_R(t)=4\left(\int |\nabla Q|^2-\int |\nabla u|^2\right)
+A_R(u(t))=-4\delta(t)+A_R(u(t)),
\end{gather}
where
\begin{equation}
 \label{E:defARradial}
A_R(u(t))= 4\int \left(\varphi''\left(\frac{x}{R}\right)-2\right)|\nabla
u|^2-\int \left(\Delta\varphi\left(\frac{x}{R}\right)-6\right) |u|^4
-\frac{1}{R^2}\int \Delta^2 \varphi\left(\frac{x}{R}\right)|u|^2.
\end{equation}

\noindent\emph{Step 1. Concavity of $y_R$.} We first claim
\begin{equation}
\label{E:y''R<0}
\exists R_0>0,\; \forall R\geq R_0,\quad y''_R(t)\leq -2\delta(t).
\end{equation}
By \eqref{E:defARradial}, we must show that there exists $R_0>0$
such that for $R\geq R_0$,
\begin{equation}
\label{E:boundARsuper}
\quad A_R(u(t))\leq 2\delta(t).
\end{equation}
The proof is close to \cite[Claim 4.3]{DuMe07} and we will only
sketch it.

Using that $e^{it}Q$ is a solution of \eqref{CP} such that the
corresponding $y_R(t)$ is constant and the corresponding $\delta(t)$
is identically zero, we get $A_R(e^{it}Q)=0$.

Recall the parameter $\delta_0$ of Section \ref{S:modulation}. First
assume that $t\in D_{\delta_1}$ (i.e., that $\delta(t)\leq
\delta_1$), where $\delta_1\in (0,\delta_0)$ will be chosen later.
By Lemma \ref{LemEstimModul1}, denoting $v=\alpha Q+h$, we get
$$
u(t)=e^{it}(Q+v(t)),\quad \|v(t)\|_{H^1}\leq C\delta(t).
$$
Noting that $\varphi''(x/R)-2=\Delta\varphi(x/R)-6= \Delta^2
\varphi(x/R)=0$ for $|x|\leq R$, we get
$$
\left|A_R(u(t))\right|=\left|A_R(Q+v)-A_R(Q)\right|\leq
C\int_{|x|\geq R}\left(Q^3|v|+|v|^4+|\nabla Q|\,|\nabla
v|+|\nabla v|^2+Q|v|+|v|^2\right)dx.
$$
In view of the exponential decay of $Q$, we obtain
$$
\left|A_R(u(t)) \right| \leq C \left(e^{-cR}\delta(t)+
\delta^2(t)+\delta^4(t)\right),
$$
which shows that there exists $R_1>0$, $\delta_1>0$ such that
\eqref{E:boundARsuper} holds for $R\geq R_1$, $t\in D_{\delta_1}$.

We now fix such a $\delta_1$ and assume $\delta(t)\geq \delta_1$.
Note that by our assumptions on $\varphi$, $\int |\nabla
u|^2(\varphi''-2)\leq 0$. It remains to bound the two other terms.
We have
\begin{equation}
 \label{E:boundARs1}
\frac{1}{R^2} \int |u|^2\Delta^2\varphi(x/R)\leq
\frac{C}{R^2}M[Q]\leq \delta_1\leq \delta(t),\quad\text{if }R\geq
R_2=\sqrt{\frac{CM[Q]}{\delta_1}}.
\end{equation}
Recall that by Strauss Lemma \cite{St77}, $u(t)$ being radial, it is
bounded and
$$
\forall x\in \RR^3\setminus\{0\},\quad |u(x,t)|\leq
\frac{C}{|x|}\|u(t)\|^{1/2}_{2}\|\nabla u(t)\|_{2}^{1/2}.
$$
Hence,
\begin{equation*}
\left|\int |u|^4(\Delta\varphi(x/R)-6)\right|
\leq C\int_{|x|\geq R}|u|^4\leq \frac{C}{R^2}\|u\|_2^3\|\nabla u\|_2
\leq \frac{C}{R^2}M[Q]^{3/2}\sqrt{\delta(t)+\|\nabla Q\|_2^2}
\end{equation*}
Using that $\delta(t)\geq \delta_1$, we get that there exists a
constant $C_{\delta_1}$, depending only on $\delta_1$ and such that
$$
\left|\int |u|^4(\Delta\varphi(x/R)-6)\right| \leq
\frac{C_{\delta_1}}{R^2}\delta(t).
$$
If $R\geq R_3=\sqrt{C_{\delta_1}}$, we get
\begin{equation}
\label{E:boundARs2}
 \left|\int |u|^4(\Delta\varphi(x/R)-6)\right|\leq \delta(t).
\end{equation}
By \eqref{E:boundARs1} and \eqref{E:boundARs2}, we get
\eqref{E:boundARsuper} for $R\geq \max\{R_2,R_3\}$ in the case
$\delta(t)\geq \delta_1$ also.

\medskip

\noindent\emph{Step 2. Proof of the finite variance of $u_0$.} Let
us fix $R\geq R_0$, where $R_0$ is given by Step 1. We first show
that for all $t$ in the domain of existence of $u$,
\begin{equation}
\label{E:y'R>0}
y_R'(t)>0.
\end{equation}
If not, using that for all $t$, $y''(t)<0$, there exists $t_1$,
$\eps>0$ such that for $t\geq t_1$, $y'_R(t)<-\eps$, which
contradicts the fact that $y$ is positive and that $u$ is globally
defined for positive time.

From the fact that $y'_R$ is positive and decreasing, we deduce that
it has a finite limit as $t$ goes to infinity. But then the integral
$\int_{0}^{+\infty} y''_R(t)dt$ is convergent, which by
\eqref{E:y''R<0} implies
$$
\int_{0}^{+\infty} \delta(s)ds<\infty.
$$
Thus, there exists a subsequence $t_n\rightarrow +\infty$ such that
$\delta(t_n)\rightarrow 0$. By Proposition \ref{P:CC}, extracting if
necessary, there exists $\theta_0\in \RR$ such that
$u(t_n)\rightarrow e^{i\theta_0}Q$ in $H^1$. By \eqref{E:y'R>0},
$y_R$ is increasing, and thus,
$$
y_R(0)=\int R^{2}\varphi(x/R)|u_0|^2\leq \int R^{2}\varphi(x/R)|Q|^2.
$$
Letting $R$ go to infinity, we get
$$
\int |x|^2|u_0|^2<\infty,
$$
which shows that we are in the finite variance case, already treated
in \S \ref{SS:supercrit-finitevariance}. \qed

\section{Convergence to $Q$ in the case
$\|\nabla u_0\|_2\|u_0\|_2<\|\nabla Q\|_2\|Q\|_2$}
\label{S:subcritical}

The main purpose of this section is to prove
\begin{prop}
 \label{P:expsub}
Consider a solution $u$ of \eqref{CP} such that
\begin{equation}
\label{equality}
M[u]=M[Q],\quad E[u]=E[Q], \quad \|\nabla u_0\|_2<\|\nabla Q\|_2.
\end{equation}
which does not scatter for positive times. Then there exists
$\theta_0 \in \cR$, $x_0 \in \cR^3$, $c, C > 0$ such that
$$
\left\|u-e^{it+i\theta_0}Q(\cdot-x_0)\right\|_{H^1}\leq Ce^{-ct}.
$$
\end{prop}
We start by proving, in \S \ref{SS:compactness} that a solution $u$
of \eqref{CP} satisfying \eqref{equality} is compact in $H^1$ up to
a translation $x(t)$ in space. In \S \ref{SS:CVmean} it is shown by
a local virial identity, that the parameter $\delta(t)=\Big|\|\nabla
u\|^2_2-\|\nabla Q\|_2^2\Big|$ converges to $0$ in mean. In \S
\ref{SS:sub_expo}, combining the results of the earlier subsections
\S \ref{SS:compactness}-\ref{SS:CVmean}, the estimates of Section
\ref{S:modulation}, and a localized virial approach with a spatial
control, we conclude the proof of Proposition \ref{P:expsub}.
Finally, \S \ref{SS:conclu_existence} is dedicated to the behavior
of the special solution $Q^-$ constructed in Proposition
\ref{R:U+U-} for negative time, concluding the proof of Theorem
\ref{T:existence}.

\subsection{Compactness properties}
\label{SS:compactness}

\begin{lemma}
 \label{L:compactness}
Let $u$ be a solution of \eqref{CP} satisfying the assumptions of
Proposition \ref{P:expsub}. Then there exists a continuous function
$x(t)$ such that
\begin{equation}
\label{E:defK} K:=\left\{u(x+x(t), t),\; t\in [0,+\infty) \right\}
\end{equation}
has a compact closure in $H^1$.
\end{lemma}
\begin{proof}[Sketch of the proof]
We only sketch the argument and refer to the proofs of \cite[Prop
4.2]{KeMe06}, \cite[Prop 5.4]{HoRo07b} and \cite[Lemma 2.8]{DuMe07}
for detailed proofs in similar contexts.

It is sufficient to show that for every time-sequence $\tau_n\geq
0$, there exists (extracting if necessary) a subsequence $x_n$ such
that $u(x+x_n,\tau_n)$ has a limit in $H^1$ (see e.g.
\cite[Appendix]{DuHoRo07P}).

By the nonradial profile decomposition
\cite[Lemma 2.1]{DuHoRo07P}, there exist families of profiles
$\psi^j\in H^1$, and of sequences $x_n^j$ and $t_n^j$ such that
\begin{gather}
 \label{E:decompo_u}
u(x,\tau_n) = \sum_{j=1}^N e^{-it_n^j\Delta} \psi^j(x-x_n^j)+W_n^M(x),
\quad \lim_{M\rightarrow+\infty} \lim_{n\rightarrow+\infty}
\|e^{it\Delta}W_n^M\|_{S(\dot{H}^{1/2})}=0,\\
 \label{E:orthogonality}
\lim_{n\rightarrow +\infty} |x_n^j-x_n^k|+|t_n^j-t_n^k|=+\infty.
\end{gather}
The crucial point is to show that there is exactly one nonzero
profile. Indeed, if for all $j$, $\psi^j=0$, then $u$ must scatter
by the local Cauchy problem theory for \eqref{CP}.

On the other hand, if at least two profiles are nonzero, then by the
Pythagorean expansions properties of the profile decomposition (see
(2.3) and (2.8) in \cite{DuHoRo07P}) there exists $\eps>0$ such that
for all $j$,
\begin{equation}
\label{E:MQpsi_j}
M[\psi^j] E[e^{-it_n^j\Delta}\psi^j]\leq M[Q]E[Q]-\eps,\quad
\|\psi^j\|_2\|\nabla \psi^j\|_2 \leq \|Q\|_2\|\nabla Q\|_2-\eps.
\end{equation}
Recall that by \cite{HoRo07b,DuHoRo07P}, a solution of \eqref{CP}
with initial condition $v_0\in H^1$ such that
$\big\|v_0\big\|_2\big\|\nabla
v_0\big\|_2<\big\|Q\big\|_2\big\|\nabla Q\big\|_2$ scatters as
$t\rightarrow \pm \infty$. By the existence of wave operators for
equation \eqref{CP}, there exists for all $j$ a function $v^j_0$ in
$H^1$ such that the corresponding solution $v^j$ of \eqref{CP}
satisfies
$$
\lim_{n\rightarrow + \infty} \left\|e^{-i t_n
^j\Delta}\psi^j-v^j(t_n^j)\right\|_{H^1}=0.
$$
Using the arguments of the proof of \cite[Prop 5.4]{HoRo07b}, one
can show, as a consequence of \eqref{E:orthogonality} and the
scattering of $v^j$, that for large $n$, the solution
$u(x,t+\tau_n)$ of \eqref{CP} is close (for positive times) to the
approximate solution $\sum_{j=1}^N v^j(x-x_n^j,t+t_n^j)$ of
\eqref{CP} (where $N$ is large). Therefore, the solution $u$ must
also scatter for positive time, which yields a contradiction,
showing that there is only one nonzero profile.

As a consequence,
$$
u(x,\tau_n)=e^{-it_n^1\Delta} \psi^1(x-x_n^1)+W_n^1(x),
\quad \lim_{n\rightarrow+\infty}\|e^{it\Delta}W_n^1\|_{S(\dot{H}^{1/2})}=0.
$$
Furthermore, $\lim_{n\rightarrow+\infty}\|W_n^1\|_{H^1}=0.$
If not, for some $\eps>0$,
$$
E\left[e^{-it_n^1\Delta} \psi^1\right]M \left[e^{-it_n^1\Delta}
\psi^1\right]\leq E[Q]M[Q]-\eps,
$$
and one can show by the preceding arguments that $u$ scatters.

It remains to show that $t_n^1$ is bounded (and thus, converges up
to extraction). If not, we may assume that $t_n^1\rightarrow+\infty$
or $t_n^1\rightarrow -\infty$. In the first case,
\begin{multline*}
\qquad \left\|e^{it\Delta}u(\tau_n)\right\|_{S\left((-\infty,0];\dot{H}^{1/2}\right)}
=\left\|e^{i(t-t_n^1)\Delta} \psi^1\right\|_{S\left((-\infty,0];
\dot{H}^{1/2}\right)}+o_n(1)\\ =\left\|e^{it\Delta}
\psi^1\right\|_{S\left((-\infty,-t_n^1];\dot{H}^{1/2}\right)}+o_n(1),\qquad
 \end{multline*}
which goes to $0$ as $n$ goes to $\infty$, showing that $u$ scatters
for positive time, a contradiction. Similarly, in the second case
$$
\left\|e^{it\Delta}u(\tau_n)\right\|_{S\left([0,+\infty);\dot{H}^{1/2}\right)}
=\left\|e^{it\Delta}
\psi^1\right\|_{S\left([-t_n^1,+\infty);\dot{H}^{1/2}\right)} +
o_n(1) \underset{n\rightarrow +\infty}{\longrightarrow}0.
$$
Thus, $u(\tau_n)$ satisfies the analogue of
\eqref{E:scattering_condition} for negative times, which shows that
$u$ scatters for negative and (by the analogue of
\eqref{E:small_scattering}),
$\|u\|_{S\left((-\infty,\tau_n];\dot{H}^{1/2}\right)}$ goes to $0$
as $n$ goes to $\infty$. Since $\tau_n\geq 0$, we get that $u=0$,
contradicting our assumptions.
\end{proof}

Let $u$ be a solution of \eqref{CP} satisfying \eqref{P:expsub}. Let
$x(t)$ be the translation parameter of Lemma \ref{L:compactness}.
Consider $\delta_0>0$ as in Section \ref{S:modulation}. The
parameters $X(t)$, $\theta(t)$, $\alpha(t)$
are defined for 
$t\in D_{\delta_0}=\{t\;|\,\delta(t)<\delta_0\}$. By
\eqref{E:decompo-u} and Lemma \ref{LemEstimModul1}, there exists a
constant $C_0>0$ such that
$$
\forall t\in D_{\delta_0}, \quad \int_{|x-X(t)|\leq 1} |\nabla
u|^2+|u|^2\geq \int_{|x|\leq 1} |\nabla
Q|^2+|Q|^2-C_0\delta(t).
$$
Taking a smaller $\delta_0$ if necessary, we can assume that the
right hand side of the preceding inequality is bounded from below by
a strictly positive constant $\eps_0$ on $D_{\delta_0}$. Thus,
$$
\forall t\in D_{\delta_0}, \quad \int_{|x+x(t)-X(t)|\leq 1} |\nabla
u(x+x(t))|^2+|u(x+x(t))|^2\geq \eps_0>0.
$$
By compactness of $\overline{K}$, it follows that
$|x(t)-X(t)|$ is bounded on $D_{\delta_0}$. As a consequence, we may
modify $x(t)$ so that $K$ defined by \eqref{E:defK} remains
precompact in $H^1$ and
\begin{equation}
\label{E:equalityX}
\forall t\in D_{\delta_0},\quad x(t)=X(t).
\end{equation}
It is classical that one may choose the function $x$ to be
continuous (see \cite[Remark 5.4]{KeMe06} and \cite[Lemma
A.3]{DuHoRo07P}). Therefore, we have shown
\begin{corol}
 \label{C:compactness}
Let $u$ be as in Proposition \ref{P:expsub}.
Then there exists a continuous function $x(t)$ such that
\eqref{E:equalityX} holds and the set $K$ defined by \eqref{E:defK}
has compact closure in $H^1$.
\end{corol}
We will also need the following:
\begin{lemma}
 \label{L:xt}
Let $u$ be as in Proposition \ref{P:expsub}, and $x(t)$ be defined
by Corollary \ref{C:compactness}. Then
\begin{equation}
 \label{E:zero_momentum}
P[u]=\im \int \overline{u}\, \nabla u \, dx=0.
\end{equation}
Furthermore,
\begin{equation}
\label{as.x}
\lim_{t\rightarrow +\infty} \frac{x(t)}{t}=0.
\end{equation}
\end{lemma}
\begin{proof}
Assume $P[u]\neq 0$ and consider, as in \cite[Prop. 4.1]{DuHoRo07P},
the Galilean transformation of $u$, $w(x,t)=e^{ix\cdot \xi_0}
e^{-it|\xi_0|^2} u(x-2\xi_0t,t)$. In order to minimize $E[w]$, we
take $\xi_0=-P[u]/M[u]$. Then $M[w]=M[u]=M[Q]$, and by the choice of
$\xi_0$, $E[w]<E[u]=E[Q]$. By the result of \cite{DuHoRo07P}, $u$
scatters in $H^1$ which contradicts our assumptions, showing
\eqref{E:zero_momentum}.

For the proof of \eqref{as.x} see \cite[Lemma 5.1]{DuHoRo07P}.
\end{proof}

\subsection{Convergence in mean}
\label{SS:CVmean}
\begin{lemma}
 \label{L:CVseq}
Let $u$ be a solution of \eqref{CP} satisfying the assumptions of
Proposition \ref{P:expsub}. Then
\begin{equation}
\label{CVint}
\lim_{T\rightarrow +\infty} \frac{1}{T} \int_0^T \delta(t) dt=0,
\end{equation}
where as in Section \ref{S:modulation}, $\delta(t)=\Big|\|\nabla
Q\|_2^2-\|\nabla u(t)\|^2_2\Big|$.
\end{lemma}
As an immediate corollary we get
\begin{corol}
 \label{C:CVseq}
Under the assumptions of Proposition \ref{P:expsub}, there exists a
sequence $t_n$ such that $t_n\rightarrow +\infty$ and
$$
\lim_{n\rightarrow+\infty} \delta(t_n)=0.
$$
\end{corol}
In the sequel we will assume, extracting if necessary, that for all
$n$, $1+t_n\leq t_{n+1}$.
\begin{proof}[Proof of Lemma \ref{L:CVseq}]

Let $\varphi$ be a $C^{\infty}$ positive radial function on $\RR^3$
such that $\varphi(x)=|x|^2$ if $|x|\leq 1$ and $\varphi(x)=0$ if
$|x|\geq 2$. Consider the localized variance
\begin{equation}
\label{defzR} y_R(t)=\int_{\RR^3} R^2 \, \varphi\left(\frac xR
\right) |u(x,t)|^2dx.
\end{equation}
Then by explicit computations and \eqref{CP},
\begin{equation}
 \label{z'R}
y_R'(t) = 2R \, \im \int \ubar \, \nabla \varphi\left(\frac
xR\right)
\cdot \nabla u,\quad |y'_R(t)|\leq CR. \\
\end{equation}
Furthermore, $y''_R(t)=\left(8\int |\nabla u|^2-6\int
|u|^4\right)+A_R(u(t))$, where
\begin{multline}
 \label{defAR}
A_R(u(t)):=4\sum_{j\neq k}\int
\frac{\partial^2\varphi}{\partial{x_j} \partial{x_k}}
\left(\frac{x}{R}\right)\frac{\partial
u}{\partial{x_j}}\frac{\partial \ubar}{\partial{x_k}}
+4 \sum_{j} \int \left(\frac{\partial^2 \varphi}{\partial x_j^2}
\left(\frac xR\right) -2\right)\left|\partial_{x_j}u\right|^2\\
 -\frac{1}{R^2} \int \Delta^2\varphi\left(\frac{x}{R}\right)|u|^2
-\int \left(\Delta \varphi\left(\frac xR\right)-6\right)|u|^4.
\end{multline}
Using as in the proof of Lemma \ref{L:supercrit-finitevariance},
that $E[u]=E[Q]$ and $M[u]=M[Q]$, we get
\begin{equation}
\label{z''R}
y''_R(t)=4\delta(t)+A_R(u(t)).
\end{equation}
Note that if $|y|\leq 1$, $(\Delta^2\varphi)(y)=0$,
$\partial_{x_j}^2\varphi(y)=2$ and $\Delta \varphi(y)=6$. Thus,
\begin{equation}
\label{bound.AR}
|A_R(u(t))|\leq C\int_{|x|\geq R} |\nabla u|^2+|u|^4+\frac{1}{R^2} |u|^2.
\end{equation}
Let $x(t)$ be as in Corollary \ref{C:compactness} and $K$ be defined by \eqref{E:defK}.
Let $\eps>0$.  By compactness of $K$, there exists $R_0(\eps)>0$ such that
\begin{equation}
 \label{defR0eps}
\forall \, t\geq 0,\quad \int_{|x-x(t)|\geq R_0(\eps)} |\nabla
u|^2+|u|^2+|u|^4\leq \eps.
\end{equation}
Furthermore, by \eqref{as.x}, there exists $t_0(\eps)\geq 0$ such that
\begin{equation}
\label{xeps}
\forall t\geq t_0(\eps), \quad |x(t)|\leq \eps t.
\end{equation}
Let
$$
T\geq t_0(\eps),\quad R:= \eps T+R_0(\eps)+1,\quad t\in [t_0(\eps),T].
$$
Let us bound the terms in \eqref{bound.AR}. Using that $|x(t)|\leq
\eps T$ and $R_0(\eps)+\eps T\leq R$, we get
\begin{multline}
 \label{bound.AR2}
\qquad\int_{|x|\geq R} |\nabla u|^2+|u|^4 +\frac{1}{R^2}|u|^2
\leq \int_{|x-x(t)|+|x(t)|\geq R} |\nabla u|^2+|u|^4+|u|^2\\
\leq \int_{|x-x(t)|\geq R_0(\eps)} |\nabla u|^2+|u|^4+|u|^2\leq\eps.
\end{multline}
By \eqref{z'R} and \eqref{z''R}, we obtain
$$
\int_{t_0(\eps)}^{T} \big[4\delta(t)+A_R(u(t))\big]dt = \int_{
t_0(\eps)}^{T} y''_R(t)dt\leq |y'_R(T)|+|y'_R(t_0(\eps))|\leq CR.
$$
Thus, by \eqref{bound.AR} and \eqref{bound.AR2}, we have, for some
constant $\widetilde{C}>0$, independent of $T$ and $\eps$,
$$
\int_{t_0(\eps)}^T \delta(t)dt
\leq C(R+T \eps)\leq \widetilde{C}\left(R_0(\eps)+1+\eps T\right).
$$
Hence,
$$
\frac{1}{T} \int_0^T \delta(t)dt \leq \frac{1}{T}\int_0^{t_0(\eps)}
\delta(t)dt+\frac{\widetilde{C}}{T}(R_0(\eps)+1)+\widetilde{C}\eps.
$$
Passing to the limit superior as $T\rightarrow+\infty$, then letting
$\eps$ tends to $0$, we get \eqref{CVint}.
\end{proof}

\subsection{Exponential convergence}
\label{SS:sub_expo}

In this section we prove Proposition \ref{P:expsub}. We refer to
\cite[Subsection 3.3]{DuMe07} and \cite[Subsection 3.3]{DuMe07Pb}
for similar arguments.

The two ingredients of the proof of Proposition \ref{P:expsub} are
the localized virial argument (Lemma \ref{L:virial}) and a precise
control of the variations of the parameter $x(t)$ (Lemma
\ref{L:parametercontrol}).

\begin{lemma}
 \label{L:virial}
Let $u$ be a solution of \eqref{CP} satisfying the assumptions of
Proposition \ref{P:expsub}, and $x(t)$ be as in Corollary
\ref{C:compactness}. Then there exists a constant $C$ such that if
$0\leq \sigma<\tau $
\begin{equation}
 \label{inegdelta}
\int_{\sigma}^{\tau} \delta(t) dt\leq C\left[1+\sup_{\sigma\leq
t\leq \tau} |x(t)|\right](\delta(\sigma)+\delta(\tau)).
\end{equation}
\end{lemma}
\begin{proof}
Consider the localized variance $y_R(t)$ defined by \eqref{defzR}.
By \eqref{z'R} and \eqref{z''R}
\begin{equation}
 \label{z'Rz''R}
y_R'(t) =2R \im \int \ubar \,\, \nabla \varphi \left(\frac xR\right)
\cdot \nabla u,\quad y''_R(t)=4\delta(t)+A_R(u(t)),
\end{equation}
where $A_R$ is defined in \eqref{defAR}.

\medskip

\noindent\emph{Step 1. Bound on $A_R$.}
We claim that if $\eps>0$, there exists a constant $R_{\eps}$ such that
\begin{equation}
 \label{E:boundAR}
\forall t\geq 0,\quad |x|\geq R_{\eps}(1+|x(t)|) \Longrightarrow
|A_R(u(t))|\leq \eps \delta(t).
\end{equation}

To prove \eqref{E:boundAR}, we distinguish two cases.

In the case when $\delta(t)$ is small, we use the estimates from
Section \ref{S:modulation}. Consider $\delta_0>0$ as in Section
\ref{S:modulation} (such that the parameters $\theta(t)$,  $X(t)$,
$\alpha(t)$ are well-defined for $t\in D_{\delta_0}$). Let
$\delta_1$ to be specified later and such that
$0<\delta_1<\delta_0$. Assume that $t\in D_{\delta_1}$. Let
$v=h+\alpha Q$, then from \eqref{E:decompo-u} and Lemma
\ref{LemEstimModul1}, we get
\begin{equation}
 \label{devuQ}
u(x,t) =e^{i(t+\theta(t))}\big[Q(x-X(t)) + v(x-X(t),t)\big]\quad
\text{and} \quad \|v(t)\|_{H^1} \leq C \delta(t).
\end{equation}
Note that if $\theta_0$ and $X_0$ are fixed, then
$e^{i\theta_0}e^{it}Q(\cdot+X_0)$ is a solution of \eqref{CP} such
that the corresponding $y_R(t)$ does not depend on $t$ and also
$\delta(t)=0$. As a consequence,
$A_R(e^{i\theta_0}e^{it}Q(\cdot+X_0))=0$ for any $R$ and $t$. By the
definition \eqref{defAR} of
$A_R$ with the change of variables $y=x-X(t)$, we obtain
\begin{align*}
|A_R(u)|
&=\left|A_R(u)-A_R\left(e^{i(t+\theta(t))}Q\left(x-X(t)\right)\right)\right|\\
&\leq C\int_{|y+X(t)|\geq R} \bigg(|\nabla Q(y)|\,|\nabla v(y)|
+|\nabla v(y)|^2\\
&\qquad \qquad \qquad +Q(y)\,|v(y)| +
|v(y)|^2 + |v(y)|^4\bigg)dy\\
&\leq C\bigg[\int_{|y+X(t)|\geq R} e^{-|y|}\left(|\nabla v(y)|+|v(y)|
+|v(y)|^3\right)\\
&\qquad \qquad+\int_{|y+X(t)|\geq R}
\left(|\nabla v(y)|^2+|v(y)|^2+|v(y)|^4\right)\bigg].
\end{align*}
By Lemma \ref{LemEstimModul1}, $\|v(t)\|_{H^1}\leq C\delta(t)$, and
hence, for some constant $C_0>0$, we get
\begin{equation}
 \label{intermAR1}
R\geq R_0+|X(t)|\Longrightarrow |A_R(u(t))| \leq
C_0\left[e^{-R_0}(\delta(t)+\delta(t)^3)+\delta(t)^2+\delta(t)^4\right].
\end{equation}
Choosing $R_0$ and $\delta_1$ such that $C_0e^{-R_0}\leq
\frac{\eps}{2}$ and $C_0(e^{-R_0}\delta_1^2+\delta_1+\delta_1^3)\leq
\frac{\eps}{2}$, we get
\begin{equation}
 \label{estim1}
R\geq R_0+|X(t)|\text{ and }\delta(t) \leq \delta_1 \Longrightarrow
|A_R(u(t))|\leq \eps \delta(t).
\end{equation}
Finally, by \eqref{E:equalityX} $x(t)=X(t)$ on $D_{\delta_0}$, which
shows that \eqref{estim1} implies \eqref{E:boundAR} for
$\delta(t)<\delta_1$.

Now assume that $\delta(t)\geq \delta_1$. Then by \eqref{defAR}, 
there exists a constant $C>0$ such that
\begin{align*}
\forall t\geq 0, \quad |A_R(u(t))|
&\leq C\int_{|x|\geq R} \left(|\nabla u(t)|^2+|u(t)|^4+|u(t)|^2\right)\\
&\leq C \int_{|x-x(t)|\geq R-|x(t)|} \left(|\nabla u(t)|^2+|u(t)|^4
+|u(t)|^2\right).
\end{align*}
By the compactness of $K$, there exists $R_{1}>0$ such that
\begin{equation}
 \label{estim2}
R\geq |x(t)|+ R_{1}\text{ and }\delta(t)\geq \delta_1\Longrightarrow
|A_R(u(t))|\leq \eps\delta_1\leq \eps\delta(t),
\end{equation}
hence, \eqref{E:boundAR} for $\delta(t)\geq \delta_1$, which
completes Step 1.

\medskip

\noindent\emph{Step 2. End of the proof.}

By \eqref{z'Rz''R} and \eqref{E:boundAR}, we get that there exists
$R_2>0$ such that
$$
R\geq R_2(1+|x(t)|)\Longrightarrow y''_R(t)\geq 2\delta(t).
$$
Let $\ds R=R_2(1+\sup_{\sigma\leq t\leq \tau}|x(t)|)$. Then
\begin{equation}
 \label{intdelta}
2 \int_{\sigma}^{\tau}\delta(t)dt\leq \int_{\sigma}^{\tau}
y''_R(t)dt=y'_R(\tau)-y'_R(\sigma)
\end{equation}
Note that if $\delta(t)<\delta_0$, then by \eqref{z'Rz''R}, \eqref{devuQ}
 and the change of variables $\xi=x-X(t)$, we get
\begin{multline*}
y_R'(t) = 2R\im \int \overline{v}(\xi) \nabla \varphi
\left(\frac {\xi+X(t)}{R}\right)\cdot \nabla Q(\xi)\\
+2 R\im \int Q(\xi) \nabla \varphi\left(\frac {\xi+X(t)}{R}\right)\cdot
\nabla v(\xi)+2R\im \int \overline{v}(\xi) \nabla \varphi\left(\frac
{\xi+X(t)}{R}\right)\cdot \nabla v(\xi),
\end{multline*}
which yields, by Lemma \ref{LemEstimModul1}, $|y'_R(t)|\leq
CR(\delta(t)+\delta(t)^2)\leq CR\delta(t)$. This inequality remains
valid if $\delta(t)\geq \delta_0$ by the straightforward estimate
$|y_R(t)|\leq CR\|\nabla u(t)\|_{2}\|u(t)\|_{2}$. In view of
\eqref{intdelta}, we get
$$
\int_{\sigma}^{\tau}\delta(t)dt\leq
CR(\delta(\sigma)+\delta(\tau))\leq C R_2 \left(1+\sup_{\sigma\leq
t\leq \tau} |x(t)|\right)(\delta(\sigma)+\delta(\tau)),
$$
which concludes the proof of Lemma \ref{L:virial}.
\end{proof}

\begin{lemma}[Control of the variations of $x(t)$]
\label{L:parametercontrol}
There exists a constant $C>0$ such that
\begin{equation}
\label{E:parametercontrol} \forall \sigma,\tau >0 \quad \text{with}
\quad \sigma+1\leq \tau, \quad |x(\tau)-x(\sigma)|\leq
C\int_{\sigma}^{\tau}\delta(t)dt.
\end{equation}
\end{lemma}

\begin{proof}
The proof is very similar to the one in \cite[Lemma 3.10]{DuMe07Pb}.
We sketch it for the sake of completeness. Let $\delta_0>0$ be as in
Section \ref{S:modulation}. Let us first show that there exist
$\delta_2 > 0$ such that
\begin{equation}
 \label{E:bounddelta}
\forall \tau \geq 0, \quad \inf_{t\in [\tau,\,\tau+2]}\delta(t)\geq
\delta_2 \quad \text{or} \quad \sup_{t\in
[\tau,\,\tau+2]}\delta(t)<\delta_0.
\end{equation}
If not, there exist two sequences $t_n,\,t'_n\geq 0$ such that
$$
\delta(t_n) \underset{n\rightarrow+\infty}{\longrightarrow 0}, \quad
\delta(t'_n)\geq \delta_0,\quad |t_n-t'_n|\leq 2.
$$
Extracting if necessary, we may assume
\begin{equation}
\label{E:CVtnt'n}
\lim_{n\rightarrow +\infty} t_n'-t_n=\tau\in [-2,2].
\end{equation}
By the compactness of $K$, $u(t_n,\cdot +x(t_n))$ converges in $H^1$
to some $v_0\in H^1$. By assumption \eqref{equality} and the fact
that $\delta(t_n)$ tends to $0$, $E[v_0]=E[Q]$, $M[v_0]=M[Q]$ and
$\|\nabla v_0\|_{2}=\|\nabla Q\|_{2}$. By Proposition \ref{P:CC},
$v_0=e^{i\theta_0}Q(\cdot-x_0)$ for some parameters $\theta_0\in
\RR$, $x_0\in \RR^3$. As a consequence, the solution of \eqref{CP}
with the initial condition $v_0$ is $e^{i(t+\theta_0)}Q(\cdot-x_0)$.
By continuity of flow and \eqref{E:CVtnt'n}, $u(t'_n,\cdot+x(t_n))$ tends to
$e^{i(\tau+\theta_0)}Q(\cdot-x_0)$ in $H^1$, which contradicts the
fact that $\delta(t_n')\geq \delta_0$, completing the proof of
\eqref{E:bounddelta}.

We now show \eqref{E:parametercontrol} with the additional condition
that $\tau\leq \sigma+ 2$. By \eqref{E:bounddelta}, we may assume
that $\sup_{t\in [\sigma,\tau]}\delta(t)< \delta_0$\, or
\,$\inf_{t\in [\sigma,\tau]} \delta(t)\geq \delta_2$. In the first
case, recalling that by the assumption \eqref{E:equalityX},
$x(t)=X(t)$ on $D_{\delta_0}$, we get
\eqref{E:parametercontrol} by time-integration of the estimate
$|X'(t)|\leq C\delta(t)$ of Lemma \ref{LemEstimModul2}. In the
second case, we have $\int_{\sigma}^{\tau}\delta(t)\geq \delta_2$
and \eqref{E:bounddelta} follows from
\begin{equation*}
\exists \, C>0,\;\forall s,t\geq 0,\quad |t-s|\leq 2\Longrightarrow
|X(t)-X(s)|\leq C,
\end{equation*}
which is a straightforward consequence of the compactness of $K$ in
$H^1$ and the continuity of the flow of equation \eqref{CP}.

To complete the proof of Lemma \ref{L:parametercontrol}, it remains
to divide $[\sigma,\tau]$ into intervals of length at least $1$ and
at most $2$ and stick together the previous inequalities to get
\eqref{LemEstimModul2} without the assumption $\tau\leq \sigma+2$.
\end{proof}

We are now ready to prove Proposition \ref{P:expsub}. Let us first
show that $x(t)$ is bounded.

Consider the sequence $\{t_n\}_{n}$ given by Corollary
\ref{C:CVseq}. Recall that $t_n$ goes to infinity, that $1+t_n\leq t_{n+1}$, and that
$\delta(t_n)$ tends to $0$. By Lemma \ref{L:virial} and Lemma
\ref{L:parametercontrol}, there exists a constant $C_0>0$ such that
$$
\forall N < n,\quad 1+t_N\leq t\leq t_n\Longrightarrow
|x(t_N)-x(t)|\leq C_0\Big(1+\sup_{t_N\leq s\leq t_n}
|x(s)|\Big)\big[\delta(t_N)+\delta(t_n)\big].
$$
Choosing $t\in [t_N+1,t_n]$ such that $|x(t)|=\sup_{t_N+1 \leq s\leq
t_n} |x(s)|$, we get
$$
\sup_{t_N+1\leq s\leq t_n} |x(s)|\leq C(N)+ C_0 \Big(1 + \sup_{t_N +
1 \leq s\leq t_n} |x(s)|\Big)\big[\delta(t_N)+\delta(t_n)\big],
$$
where $\ds C(N)=|x(t_N)|+C_0\sup_{t\in [t_N,t_N+1]}|x(t)|$ (we
assumed $\delta(t_N)+\delta(t_n)\leq 1$). Fix $N$ large enough such
that $C_0\delta(t_N)\leq \frac 12$. Then as soon as $t_n\geq t_N+1$,
$$
\frac{1}{2}\sup_{t_N+1\leq s\leq t_n} |x(s)|\leq C(N)+ \frac
12+C_0\Big(1+\sup_{t_N+1\leq s\leq t_n} |x(s)|\Big)\delta(t_n).
$$
Letting $n$ tend to infinity and using again that $\delta(t_n)$
tends to $0$, we get that $|x(t)|$ is bounded on $[t_N+1,+\infty)$,
and thus, by continuity, on $[0,+\infty)$.

We will now show that
\begin{equation}
 \label{E:expint}
\exists \, c,C>0,\;\forall \sigma\geq 0,\quad
\int_{\sigma}^{\infty}\delta(t)dt\leq Ce^{-c\sigma},
\end{equation}
which will yield, together with Lemma \ref{L:expo}, the conclusion
of Proposition \ref{P:expsub}.

By Lemma \ref{L:virial} and the boundedness of $x(t)$,
\begin{equation*}
\forall \, \sigma, \tau > 0 ~~\text{ such that }~~0\leq \sigma\leq
\tau,\quad \int_{\sigma}^{\tau} \delta(t)dt\leq
C\left(\delta(\sigma)+\delta(\tau)\right).
\end{equation*}
Fix $\sigma$ and take $\tau=t_n$, where the sequence $(t_n)_n$,
given by Corollary \ref{C:CVseq} is such that $\lim_n
\delta(t_n)=0$. Letting $n$ tend to $\infty$, we get that
$\int_0^{+\infty}\delta(t)dt$ is finite and for $\sigma\geq 0$,
$\int_{\sigma}^{+\infty}\delta(t)dt\leq C\delta(\sigma)$. Gronwall's
Lemma yields \eqref{E:expint}, concluding the proof of Proposition
\ref{P:expsub}.\qed

\subsection{Scattering of $Q^-$ for negative times.}
\label{SS:conclu_existence}

In this paragraph, we conclude the proof of Theorem
\ref{T:existence} by showing by contradiction that the special
solution $Q^-$ constructed in Proposition \ref{P:UA} and Remark
\ref{R:U+U-} scatters as $t\rightarrow -\infty$.

If not, applying the arguments of \S
\ref{SS:compactness}-\ref{SS:sub_expo} to the solutions $Q^-$ and
$t\rightarrow \overline{Q^{-}}(x,-t)$ of \eqref{CP}, we get that
there exists a parameter $x(t)$, defined for $t\in \RR$ and such
that $\widetilde{K}=\left\{Q^-(\cdot+x(t),t),\;t\in\RR\right\}$ has
compact closure in $H^1$. By the argument at the end of \S
\ref{SS:sub_expo}, $x(t)$ is bounded and $\delta(t)$ tends to $0$ as
$t\rightarrow \pm \infty$. A simple adjustment of Lemma
\ref{L:virial} gives
\begin{equation*}
-\infty<\sigma\leq \tau<+\infty\Longrightarrow
\int_{\sigma}^{\tau}\delta(t)dt\leq C\left[1 + \sup_{\sigma\leq
t\leq \tau} |x(t)|\right](\delta(\sigma)+\delta(\tau))\leq
C(\delta(\sigma)+\delta(\tau)).
\end{equation*}
Letting $\sigma$ go to $-\infty$ and $\tau$ to $+\infty$, we get
$\int_{\RR} \delta(t)dt=0$, thus, $\delta(t)=0$ for all $t$,
contradicting the assumption $\|\nabla u_0\|_2<\|\nabla Q\|_2$.

\section{Uniqueness}
\label{S:uniqueness} In this section, to conclude the proof of
Theorem \ref{T:classification}, we show the following uniqueness
statement:
\begin{prop}
 \label{P:uniqueness}
Let $u$ be a solution of \eqref{CP}, defined on $[t_0,+\infty)$,
such that $E[u]=E[Q]$, $M[u]=M[Q]$ and
\begin{equation}
 \label{E:expuniqueness}
\exists \, c,C>0 : \; \quad \|u-e^{it} Q\|_{H^1}\leq Ce^{-ct} \quad
\forall t\geq t_0.
\end{equation}
Then there exists $A\in \RR$ such that $u=U^A$, where $U^A$ is the
solution of \eqref{CP} defined in Proposition \ref{P:UA}.
\end{prop}

The proof of Proposition \ref{P:uniqueness} relies on a careful
analysis of solutions of the linearized equation (equation
\eqref{eqvg} below), that decay exponentially as $t$ tends to
$+\infty$. This analysis, carried out in \S \ref{SS:expo}, relies on
the spectral properties of $\LLL$ described in \S
\ref{SS:linearized}. In \S \ref{SS:uniqueness} we finish the proof
of Proposition \ref{P:uniqueness}, and in \S
\ref{SS:end_classification} we gather the results of Sections
\ref{S:supercritical}, \ref{S:subcritical} and \ref{S:uniqueness} to
show Theorem \ref{T:classification}.

\subsection{Exponentially small solutions of the linearized equation}
\label{SS:expo}

Recall the notation of Section \ref{S:spectralsolution}, in
particular the operator $\LLL$ and its eigenvalues and
eigenfunctions. Consider
$$
v\in C^{0} \left([t_0,+\infty),H^1\right) \quad \text{and} \quad
g\in C^{0}\left([t_0,+\infty),L^{2}\right)
$$
such that
\begin{gather}
 \label{eqvg}
\partial_t v+\LLL v=g,\quad (x,t)\in \RR^3\times (t_0,+\infty), \\
 \label{decexp}
\|v(t)\|_{H^1}\leq Ce^{-\gamma_1t},\quad \|g(t)\|_{2}\leq
Ce^{-\gamma_2 t}, \quad t\geq t_0,
\end{gather}
where
\begin{equation*}
0<\gamma_1<\gamma_2.
\end{equation*}

For any $\gamma \in \cR$, denote by $\gamma^-$ a positive number
arbitrary close to $\gamma$ and such that $0<\gamma^-<\gamma$.

We now prove the following self-improving estimate.
\begin{lemma}
 \label{L:Linearized}
Under the above assumptions,
\begin{enumerate}
\item
 \label{I:1}
if $\gamma_2\leq e_0$, then $\|v(t)\|_{H^1}\leq Ce^{-\gamma_2^- t}$,

\item
 \label{I:2}
if $\gamma_2>e_0$, then there exists $A\in \RR$ such that
$v(t)=Ae^{-e_0t}\YYY_++w(t)$ with $\|w(t)\|_{H^1}\leq
Ce^{-\gamma_2^- t}$.
\end{enumerate}
\end{lemma}

\begin{proof}
In this proof we work with the real $L^2$-scalar product, denoted by
$(\cdot,\cdot)$,
$$
(u,v)=\re \int u\,\overline{v}=\int \re u \re v+\int \im u\im v.
$$
We first normalize the eigenfunctions of $\LLL$. Denote
$$
Q_0 :=\frac{iQ}{\left\|Q\right\|_{2}},\quad
Q_j:=\frac{\partial_{x_j}Q}{\left\|\partial_{x_j}Q\right\|_{2}}.
$$
From \S \ref{SS:linearized} recall the quadratic form on $H^1$,
$\Phi$, and its associated bilinear form $B$. From
\eqref{E:orthobilinear} we have
$$
\forall j=0, \ldots 3,\quad \forall h\in H^1,\quad B(Q_j,h)=0,\quad
\|Q_j\|_{2}=1.
$$
Recall that $\Phi(\YYY_+)=\Phi(\YYY_-)=0$ and $B(\YYY_+,\YYY_-)\neq
0$. Normalize the eigenfunctions $\YYY_+$, $\YYY_-$ such that
$B(\YYY_+,\YYY_-)=1$. We have
$$
h \in G_{\bot}'\iff \forall j,\quad (Q_j,h)=0 \text{ and
}B(\YYY_+,h)=B(\YYY_-,h)=0.
$$
Indeed, if $h_1=\re h$, $h_2=\im h$, in view of \eqref{bilinear} and Remark \ref{R:Y1Y2},
$$
B(\YYY_+, h) = \frac{e_0}{2} \Big[\left(\YYY_2,h_1\right) - \left(
\YYY_1, h_2\right) \Big],\quad B(\YYY_ -, h) = \frac{e_0}{2} \Big[
\left(\YYY_2,h_1\right)+\left(\YYY_1,h_2\right)\Big],
$$
which shows that the orthogonality condition \eqref{orthoef} is
equivalent to the condition $B(\YYY_+,h)=B(\YYY_-,h)=0$.

Next, write
\begin{gather}
 \label{decompov}
v(t) =\alpha_+(t)\YYY_+ +\alpha_-(t)\YYY_- +\sum_{j=0}^3
\beta_j(t)Q_j + v_{\bot}(t),\quad v_{\bot}\in G_{\bot}',\\
 \label{expralpha}
\text{where }\alpha_{+}(t)=B(v(t),\YYY_-),\quad \alpha_{-}(t)=B(v(t),\YYY_+)\\
 \label{exprbeta}
\forall j\in \{0,\ldots,3\},\quad \beta_j(t)
=(v(t),Q_j)-\alpha_+(t)(\YYY_+,Q_j)-\alpha_-(t)(\YYY_-,Q_j).
\end{gather}
By the radiality of $\YYY_{\pm}$ and $Q$, we have
$(\YYY_+,Q_j)=(\YYY_-,Q_j)=0$ for $j=1,2,3$, but we will not need
this property in the sequel.
\medskip

\noindent\emph{Step 1. Differential equations on the coefficients.} Let us show
\begin{gather}
\label{diff1}
\left|\alpha_-'(t)-e_0\alpha_-(t)\right|\leq Ce^{-\gamma_2t},
\quad \left|\alpha_+'(t)+e_0\alpha_+(t)\right|\leq Ce^{-\gamma_2t}.\\
\label{diff2}
\forall j\in \{0,\ldots,3\},\quad |\beta'_j(t)|
\leq C\left(\|v_{\bot}(t)\|_{2}+e^{-\gamma_2t}\right)\\
\label{diff3}
\left|\frac{d}{dt}\Phi(v(t))\right|\leq e^{-(\gamma_1+\gamma_2)t}
\end{gather}
First note that $\LLL$ is antisymmetric for the bilinear form $B$.
Indeed, for $g,h\in H^2$ 
$$
B(g,\LLL h)=\frac 12(L_+ g_1,-L_- h_2)+\frac 12(L_-g_2,L_+ h_1) = -
B(\LLL g,h).
$$
By \eqref{eqvg} and \eqref{expralpha}, we have
\begin{equation*}
\alpha'_-(t)=B(\partial_t v,\YYY_+)=B(-\LLL
v+g,\YYY_+)=B(v,\LLL\YYY_+) +B(g,\YYY_+)=e_0\alpha_-(t)+B(g,\YYY_+).
\end{equation*}
In view of assumption \eqref{decexp} on $g$, we get the inequality
on $\alpha_-(t)$ in \eqref{diff1}. The inequality on $\alpha_+(t)$
follows from the same argument.

By \eqref{exprbeta}, we obtain
\begin{align*}
\beta_j'
&=\big(\partial_t v-\alpha_+'\YYY_+-\alpha_-'\YYY_-,Q_j\big)
=\big(-\LLL v-\alpha_+'\YYY_+-\alpha_-'\YYY_-,Q_j\big)+ \big(g,Q_j\big)\\
&=\big(-\alpha_+\LLL \YYY_+-\alpha_-\LLL\YYY_--\alpha_+'\YYY_+-\alpha_-'\YYY_-,
Q_j\big)-\big(\LLL v_{\bot},Q_j\big)+\big(g,Q_j\big)
\end{align*}
Applying \eqref{diff1}, the first term above is estimated as
$$
\big| \big(-\alpha_+\LLL \YYY_+ -\alpha_-\LLL\YYY_- -\alpha_+'\YYY_+
-\alpha_-'\YYY_-, Q_j\big) \big| \leq Ce^{-\gamma_2 t}.
$$
The assumption \eqref{decexp} implies $|\big(g,Q_j\big)|\leq
Ce^{-\gamma_2 t}$. Furthermore, $(\LLL v_{\bot},Q_j) =
\left(v_{\bot},
\LLL^*Q_j\right)$, where $\LLL^*:= \begin{pmatrix} 0& L_+\\
-L_-
&0\end{pmatrix}$ is the $L^2$-adjoint of $\LLL$, which shows the
estimate $\left|\left(\LLL v_{\bot},Q_j\right)\right|\leq
C\left\|v_{\bot}\right\|_{2}$, completing the proof of
\eqref{diff2}.

It remains to prove \eqref{diff3}. We have
\begin{equation*}
\frac{d}{dt}\Phi(v(t))= 2B \big(\partial_t v(t), v(t) \big) = -2B
\big( \LLL v,v\big)+2B(g,v).
\end{equation*}
As $B(\LLL v,v)=-B(\LLL v,v)$, we get that $B(\LLL v,v)=0$, which
yields \eqref{diff3}, using again the assumption \eqref{decexp} on
$g$, and hence, completing Step 1.
\medskip

\noindent\emph{Step 2.}
Let us show
\begin{alignat}{2}
 \label{gamma2<e0}
|\alpha_+(t)|&\leq Ce^{-\gamma_2^-t}&\quad&\text{ if }\gamma_2\leq e_0,\\
 \label{gamma2>e0}
\exists A\in \RR,\quad |\alpha_+(t)-Ae^{-e_0t}|
&\leq C e^{-\gamma_2 t}&\quad&\text{ if }\gamma_2> e_0.
\end{alignat}
Indeed, by the second inequality in \eqref{diff1}, we obtain
\begin{equation}
 \label{diffstep2}
\left|\frac{d}{dt}(e^{e_0t}\alpha_+(t))\right|\leq Ce^{(e_0-\gamma_2)t}.
\end{equation}
First assume that $\gamma_2\leq e_0$. Then by \eqref{diffstep2}, for $t\geq t_0$,
$$
\left|e^{e_0 t}\alpha_+(t)\right| \leq \begin{cases}
e^{e_0t_0}\alpha_+(t_0)+Ce^{(e_0-\gamma_2)t}&\text{ if } e_0>\gamma_2\\
e^{e_0t_0}\alpha_+(t_0)+C(t-t_0)&\text{ if } e_0=\gamma_2,
                                          \end{cases},
$$
which gives \eqref{gamma2<e0}.

Now assume $\gamma_2>e_0$. Then $\int_{t_0}^{\infty}
e^{(e_0-\gamma_2)t}dt<\infty$. By \eqref{diffstep2}, we get that
$e^{e_0 t}\alpha_+(t)$ has a limit $A$ as $t\rightarrow \infty$ and
$$
\left|e^{e_0t}\alpha_+(t) -A\right|\leq Ce^{(e_0-\gamma_2)t},
$$
implying \eqref{gamma2>e0}.
\medskip

\noindent\emph{Step 3. Conclusion of the proof in a reduced case.}
Let us conclude the proof when $\gamma_2\leq e_0$, or when
$\gamma_2>e_0$ and $A=0$. In both cases we have, in view of
\eqref{diff1}, \eqref{gamma2<e0} and \eqref{gamma2>e0},
\begin{equation}
 \label{inegalpha+}
\forall t\geq t_0,\quad |\alpha_+(t)|+|\alpha_+'(t)|\leq Ce^{-\gamma_2^-t}.
\end{equation}
By the first inequality in \eqref{diff1},
$$
\left|\frac{d}{dt}(e^{-e_0t}\alpha_-(t))\right|\leq Ce^{-(e_0+\gamma_2)t}.
$$
Integrating between $t$ and $+\infty$, we get $|\alpha_-(t)| \leq
Ce^{-\gamma_2t}$, and by \eqref{diff1}, it follows that
\begin{equation}
 \label{inegalpha-}
\forall t\geq t_0,\quad |\alpha_-(t)|+|\alpha_-'(t)|\leq Ce^{-\gamma_2t}.
\end{equation}
By the decomposition \eqref{decompov} of $v$, and recalling that as
$v_{\bot}\in G_{\bot}'$, we have
$B(\YYY_+,v_{\bot})=B(\YYY_-,v_{\bot})=0$, and that
$B(\YYY_+,\YYY_-)=1$, $B(\YYY_+,\YYY_+)=B(\YYY_-,\YYY_-)=0$, we get
$$
\Phi(v)=B(v,v)=B(v_{\bot},v_{\bot})+2\alpha_+\alpha_-.
$$
By \eqref{diff3}, \eqref{inegalpha+} and \eqref{inegalpha-},
$$
\left|\frac{d}{dt}B (v_{\bot},v_{\bot})\right| \leq C \big( e^{-2
\gamma_2^-t} + e^{-(\gamma_1+\gamma_2)t}\big) \leq C e^{-(\gamma_1 +
\gamma_2) t}.
$$
Noting that $B (v_{\bot},v_{\bot})\rightarrow 0$ as $t\rightarrow
+\infty$, and integrating the previous inequality between $t$ and
$+\infty$, we get $|B (v_{\bot},v_{\bot})|\leq
Ce^{-(\gamma_1+\gamma_2)t}$. By Proposition \ref{P:spectralL}, we
obtain
\begin{equation}
 \label{inegvbot}
\forall t\geq t_0,\quad \left\|v_{\bot}(t)\right\|_{H^1} \leq
Ce^{-\left(\frac{\gamma_1+\gamma_2}{2}\right)t}.
\end{equation}
Let $j\in \{0,\ldots,3\}$. By \eqref{diff2}, and \eqref{inegvbot},
$$
|\beta'_j(t)|\leq
C\left(e^{-\left(\frac{\gamma_1+\gamma_2}{2}\right)t} + e^{
-\gamma_2 t} \right)\leq C e^{-\left(
\frac{\gamma_1+\gamma_2}{2}\right)t}.
$$
Integrating again between $t$ and $+\infty$, we get
\begin{equation}
 \label{inegbeta}
\forall \, t\geq t_0,\quad \left|\beta_j(t)\right|\leq
Ce^{-\left(\frac{\gamma_1+\gamma_2}{2}\right)t}.
\end{equation}
In view of the decomposition \eqref{decompov} of $v$, the
inequalities \eqref{inegalpha+}, \eqref{inegalpha-},
\eqref{inegvbot} and \eqref{inegbeta} imply
$$
\forall \, t \geq t_0,\quad \|v(t)\|_{H^1}\leq
Ce^{-\left(\frac{\gamma_1+\gamma_2}{2}\right)t}.
$$
Thus, $v$ and $g$ satisfy the assumptions \eqref{decexp}, with
$\gamma_1$ replaced by $\gamma_1'=\frac{\gamma_1+\gamma_2}{2}$. An
iteration argument yields
\begin{equation}
 \label{finalinegv}
\|v(t)\|_{H^1}\leq Ce^{-\gamma_2^- t},
\end{equation}
which concludes the proof when $\gamma_2\leq e_0$ or $A=0$.
\medskip

\noindent\emph{Step 4. Conclusion of the proof in the case
$\gamma_2>e_0$, $A\neq 0$.}
Note that if $\gamma_1>e_0$, we must
have $A=0$, so that we may assume $\gamma_1\leq e_0$. Let
$$
\tv(t)=v(t)-Ae^{-e_0t}\YYY_+.
$$
Then
\begin{equation*}
\partial_t\tv(t)+\LLL \tv(t)=g(t),\quad \|\tv(t)\|_{H^1}\leq Ce^{-\gamma_1t},
\end{equation*}
and by \eqref{gamma2>e0},
\begin{equation*}
\lim_{t\rightarrow+\infty} e^{e_0t}\tilde{\alpha}_+(t)=0,
\end{equation*}
where $\tilde{\alpha}_+(t)=B(\tv(t),\YYY_-)$ is the coefficient of
$\YYY_+$ in the decomposition of $\tv(t)$ analogous to
\eqref{decompov}. Thus, $\tv(t)$ and $g$ satisfy all the assumptions
of Step 3. Hence,
$$
\left\|v(t)-Ae^{-e_0 t}\YYY_+\right\|_{H^1}\leq Ce^{-\gamma_2^-t},
$$
which concludes the proof of Lemma \ref{L:Linearized} in this case also.
\end{proof}

\subsection{Uniqueness}
\label{SS:uniqueness} Let us prove Proposition \ref{P:uniqueness}.
Let $u$ satisfy the hypothesis and write $u=e^{it}(Q+h).$

\noindent\emph{Step 1. Improvement of the decay at infinity.} We start with
showing that if $e_0^-$ is any positive number such that
$e_0^-<e_0$,
\begin{equation}
\label{E:decaye0-}
\forall t\geq t_0,\quad \|h(t)\|_{H^1}\leq Ce^{-e_0^-t}.
\end{equation}
Indeed, we have $\partial_t h+\LLL h=R(h)$, where the remainder term
$R(h)$, defined in \eqref{defR}, is a sum of quadratic and cubic
terms in $h$. By the assumption \eqref{E:expuniqueness} and Sobolev
embeddings, $\|h(t)\|_{p} \leq Ce^{-ct}$  for every $p\in [2,6]$,
which yields the bound $\|R(h)\|_2\leq Ce^{-2ct}$.  Thus, $h$
satisfies the assumptions of Lemma \ref{L:Linearized} with $g=R(h)$,
$\gamma_1=c$, $\gamma_2=2c$. If $2c>e_0$, the proof is complete. If
not, we get $\|h(t)\|_{H^1}\leq Ce^{-2c^-t}$, and the result follows
from an iteration argument.
\medskip

\noindent\emph{Step 2.} Consider the special solutions $U^A$
constructed in Proposition \ref{P:UA}, and write
$U^A=e^{it}(Q+h^A)$. Let us show that there exists $A\in \RR$ such
that for all $\gamma>0$,
\begin{equation}
 \label{E:nearhA}
\exists \, C>0, \;\forall \, t\geq t_0,\quad
\|h(t)-h^A(t)\|_{H^1}\leq Ce^{-\gamma t}.
\end{equation}
According to Step 1, $h$ fulfills the assumptions of Lemma
\ref{L:Linearized} with $\gamma_1=e_0^-$, $\gamma_2=2e_0^-$. Thus,
there exists $A\in \RR$ such that
\begin{equation}
 \label{E:asymptotich}
\left\|h(t)-Ae^{-e_0t}\YYY_+\right\|_{H^1}\leq Ce^{-2e_0^-t}.
\end{equation}
By the asymptotic development of $h^A$ obtained in Section 3,
$$
\left\|h^A(t)-Ae^{-e_0t}\YYY_+\right\|_{H^1}\leq Ce^{-2e_0 t}.
$$
Thus, \eqref{E:asymptotich} yields \eqref{E:nearhA} for any $\gamma
< 2 e_0$. We next show that if \eqref{E:nearhA} holds for some
$\gamma>e_0$, it also holds for $\gamma'=\gamma+\frac 12 e_0$. Note
that $h-h^A$ is a solution to the equation
\begin{equation*}
\partial_t (h-h^A)+\LLL (h-h^A)=R(h)-R(h^A).
\end{equation*}
By the explicit expression of $R$, and Sobolev inequalities, we get
$$
\left\|R(h)-R(h^A)\right\|_{2} \leq C
\|h-h^A\|_{H^1}\Big(\|h\|_{H^1}+\|h^A\|_{H^1}+\|h\|_{H^1}^2+\|h^A\|_{H^1}^2\Big).
$$
If \eqref{E:nearhA} holds for some $\gamma>e_0$, then
$$\left\|R(h)-R(h^A)\right\|_{2}\leq Ce^{-(e_0+\gamma)t},$$
which shows that $h-h^A$ fulfills the assumptions of Lemma
\ref{L:Linearized} with $\gamma_1=\gamma$, $\gamma_2=\gamma+e_0$,
yielding \eqref{E:asymptotich} with $\gamma+\frac 12 e_0$ instead of
$\gamma$. Step 2 is complete.

\medskip

\noindent\emph{Step 3. Uniqueness argument.}

We are now ready to finish the proof of Proposition
\ref{P:uniqueness}. Let $v:=h-h^A$. We must show that $v=0$. We will
use that $v$ is a solution to the following Schr\"odinger equation
\begin{equation}
\label{E:uniquenessv} i\partial_t v+\Delta v-v=-2Q^2 v-i Q^2\overline{v}+M,
\end{equation}
where $M(t)=iR(h(t))-iR(h^A(t))$. By H\"older's inequality
and the decay of $h$ and $h^A$ at infinity, there exists a constant
$C_1>0$ such that
\begin{equation}
\label{E:estimateM}
\forall t\geq t_0,\quad \|M(t)\|_{6/5}\leq C_1e^{-e_0t}\|v(t)\|_{2}.
\end{equation}
Let $t_1\geq t_0$, $\tau>0$ and $I=(t_1,t_1+\tau)$. By Strichartz
estimates, there exists $K>0$ such that
\begin{equation*}
\|v\|_{L^{\infty}(I;L^2)}\leq K\left\{\|v(t_1+\tau)\|_{2} +\|Q^2
v\|_{L^1(I,L^2)}+\|M(t)\|_{L^2(I;L^{6/5})}\right\}.
\end{equation*}
Integrating in time on $I$ the square of \eqref{E:estimateM}, we get
$\|M(t)\|_{L^2(I,L^{6/5})}\leq
\frac{C_1}{\sqrt{2e_0}}e^{-e_0t_1}\|v\|_{L^{\infty}(I,L^2)}$.
Furthermore, $\|Q^2 v\|_{L^1(I,L^2)}\leq
\tau\|Q^2\|_{\infty}\|v\|_{L^{\infty}(I,L^2)}$. Hence,
\begin{equation*}
\|v\|_{L^{\infty}(I;L^2)}\leq K\left\{\|v(t_1+\tau)\|_{2}
+\tau\|Q^2\|_{\infty}\|v\|_{L^{\infty}(I,L^2)}
+\frac{C_1}{\sqrt{2e_0}}e^{-e_0t_1}\|v\|_{L^{\infty}(I,L^2)}\right\}.
\end{equation*}
Let $\tau=\frac{1}{3K\|Q^2\|_{\infty}}$, choose $T\geq t_0$ such
that $\frac{C_1}{\sqrt{2e_0}}e^{-e_0T}\leq \frac{1}{3K}$. Then for
$t_1\geq T$,
\begin{equation*}
\|v(t_1)\|_{2}\leq \|v\|_{L^{\infty}(I;L^2)}\leq 3K\|v(t_1+\tau)\|_{2}.
\end{equation*}
By induction we get
$$
(3K)^n\|v(T)\|_{2}\leq \|v(T+n\tau))\|_{2},
$$
which contradicts \eqref{E:nearhA} if $\gamma$ is large enough,
unless $v(T)=0$. Thus, $h(T)=h^A(T)$     and by uniqueness in
\eqref{CP}, $h=h^A$, and thus, $u=U^A$, concluding the proof of
Proposition \ref{P:uniqueness}. \qed

\subsection{Proof of the classification result}
\label{SS:end_classification}

In this subsection we prove Theorem \ref{T:classification}.

We first show that if $A\neq 0$, $U^A$ is equal to $Q^{+}$ (if
$A>0$) or $Q^-$ (if $A<0$) up to a translation in time and a
multiplication by a complex number of modulus $1$. Indeed, by
\eqref{E:smallUa} and the definition of $Q^{\pm}$ in Remark
\ref{R:U+U-}, we have
\begin{equation}
 \label{E:asyU+-}
Q^{\pm}(t) = e^{it}Q \pm e^{-e_0 t_0}\, e^{(i-e_0)t}
\YYY_++O\left(e^{-2e_0t}\right)\text{ in } H^1.
\end{equation}
Fix $A>0$ (the proof is similar when $A<0$). Let
$t_1=-t_0-\frac{1}{e_0} \log A$, so that $e^{-e_0(t_0+t_1)}=A$. By
\eqref{E:smallUa} and \eqref{E:asyU+-}, we obtain
\begin{equation}
\label{E:devtU+} e^{-it_1}Q^{+}(t+t_1)=e^{it}Q
+e^{-e_0(t_0+t_1)}e^{-e_0t}e^{it}\YYY_++
O\left(e^{-2e_0t}\right)=U^{A}+O(e^{-2e_0t})\text{ in }H^1.
\end{equation}
As a consequence $e^{-it_1}Q^+(t+t_1)-e^{it}Q$ tends to $0$
exponentially in $H^1$ as $t\rightarrow +\infty$. By Proposition
\ref{P:uniqueness}, there exists $\widetilde{A}$ such that
$e^{-it_1}Q^+(t+t_1)=U^{\widetilde{A}}$. By \eqref{E:devtU+} we have
$\widetilde{A}=A$, which shows that $U^A=e^{-it_1}Q^{+}(t+t_1)$.

Let $u$ be a solution of \eqref{CP} satisfying the assumptions of
Theorem \ref{T:classification}. Then $M[u]E[u]=M[Q]E[Q]$. Rescaling
$u$ we may assume
$$
M[u]=M[Q],\quad E[u]=E[Q].
$$

If $\|\nabla u_0\|_2=\|\nabla Q\|_2$ (case \eqref{case_critical}),
then by the variational characterization of $Q$ (see \S \ref{SS:Q})
$u_0(x)=e^{i\theta_0} Q(x-x_0)$ for some parameters $\theta_0$,
$x_0$, and thus, by uniqueness of the Cauchy problem \eqref{CP},
$u(x,t)=e^{i\theta_0+it}Q(x-x_0)$.  Thus, $u$ is equal to $e^{it}Q$
up to the symmetries of the equation, yielding case
\eqref{case_critical}.

Assume next $\|\nabla u_0\|_2< \|\nabla Q\|_2$ (case
\eqref{case_sub}). By assumption, $u$ does not scatter for both
positive and negative times. Replacing $u(x,t)$ by
$\overline{u}(x,-t)$ if necessary, we may assume that $u$ does not
scatter for positive times. By Proposition \ref{P:expsub}, there
exists $\theta_0 \in \cR$, $x_0 \in \cR^3$, and $c, C > 0$ such that
$$
\left\|u(t)-e^{it+i\theta_0}Q(\cdot-x_0)\right\|_{H^1}\leq
Ce^{-c\,t}, ~~ t>0.
$$
Hence, $v(x,t)=e^{-i\theta_0}u(x+x_0,t)$ satisfies the assumptions
of Proposition \ref{E:expuniqueness}, which shows that $v=U^A$ for
some $A$. As $\|\nabla u_0\|_2<\|\nabla Q\|_2$, the parameter $A$
must be negative proving that $v$ (and thus $u$) is equal to $Q^-$
up to the symmetries of the equation. Therefore, case
\eqref{case_sub} of the theorem follows.

The proof of case \eqref{case_super}, combining Proposition
\ref{P:supercrit} and Proposition \ref{E:expuniqueness}, is similar
to the proof of case \eqref{case_sub} and left to the reader. \qed

\appendix

\section{Coercivity properties of the quadratic form}

 \label{A:coercivity}

This appendix is dedicated to the proof of the results of \S
\ref{SS:linearized}. Before proving Proposition \ref{P:Coercivity},
we show \eqref{BoundPhi}. Consider $h\in H^1$ and assume
$E[Q+h]=E[Q]$ and $M[Q+h]=M[Q]$. Expanding $E[Q+h]$ in terms of $Q$
and $h$, we get
$$
E[Q+h]=E[Q]+\int \nabla Q\cdot \nabla h_1-\int Q^3 h_1+\frac 1 2 \int
|\nabla h|^2 -\frac{1}{2}\int Q^2(3h_1^2+h_2^2)-\int Q |h|^2 h_1 -
\frac14 \int |h|^4.
$$
Since $E[Q+h]=E[Q]$ and $\int \nabla Q\cdot \nabla h_1-\int Q^3
h_1=-\int (\Delta Q+Q^3) h_1=-\int Qh_1$ by \eqref{E:Q}, we obtain
$$
0 = -\int Q h_1+\frac 12 \int |\nabla h|^2-\frac 12 \int Q^2
(3h_1^2+h_2^2)-\int Q |h|^2 h_1 - \frac14 \int |h|^4.
$$
Furthermore, $M[Q+h]=M[Q]$ implies $2\int Q h_1+\int |h|^2=0$,
yielding
\begin{equation}
 \label{defPhi}
\Phi(h)=\int Q |h|^2 h_1+\frac14 \int |h|^4,
\end{equation}
which gives \eqref{BoundPhi}.

The remainder of the Appendix is dedicated to the proof of
Proposition \ref{P:Coercivity}.

\subsection{Coercivity of $\Phi$ on $G_{\bot}$.}

Let us prove \eqref{coercivity} when $h\in G_{\bot}$ (see
\cite{We85,We86b}, \cite[ex. B11-B14]{Ta06BO} for similar proofs for
mass-subcritical NLS). We divide the proof into two steps.

\noindent\emph{Step 1. Nonnegativity.} We show, as a consequence of
Gagliardo-Nirenberg inequality \eqref{E:GN}, that if $h\in H^1$
satisfies \eqref{ortho2}, then
\begin{equation}
\label{positivity}
\Phi(h)\geq 0.
\end{equation}
For $u\in H^1$, let
\begin{equation}
 \label{E:I}
I(u):=\frac{\|\nabla u\|^3_2 \|u\|_2}{\|\nabla Q\|^3_2
\|Q\|_2}-\frac{\|u\|^4_4}{\|Q\|^4_4}.
\end{equation}
By \eqref{E:GN}, $I(u)\geq 0$. Take $h\in H^1$, $\alpha\in \RR$ and
compute the expansion of $I(Q+\alpha h)$ in $\alpha$ of order $2$.
By \eqref{ortho2}, we have $\int \nabla Q\cdot\nabla h_1=0$, and
thus,
\begin{align*}
\left(\int |\nabla (Q+\alpha h)|^2\right)^{\frac 32}
&=\left(\int |\nabla Q|^2\right)^{3/2}\Bigg(1+\frac 32\alpha^2\frac{\int |\nabla h|^2}{\int |\nabla Q|^2} +O\big(\alpha^4\big)\Bigg).
\intertext{Furthermore,}
\left(\int |Q+\alpha h|^2\right)^{1/2}
&=\left(\int Q^2\right)^{1/2}\left(1+\alpha \frac{\int Q h_1}{\int Q^2}+\frac 12\alpha^2\frac{\int |h|^2}{\int Q^2}-\frac 12\alpha^2\frac{(\int Qh_1)^2}{\left(\int Q^2\right)^2}+O\big(\alpha^3\big) \right)
\intertext{and}
\int |Q+\alpha h|^4
&=\left(\int Q^4\right)\left(1+4\alpha \frac{\int Q^3 h_1}{\int Q^4} +\alpha^2\frac{\int 6Q^2 h_1^2+\int 2Q^2 h_2^2}{\int Q^4}\right)+O\Big(\alpha^3\Big).
\end{align*}
Substituting above quantities into \eqref{E:I}, we obtain
\begin{multline*}
I(Q+\alpha h) = \left( 1+ \frac{3}{2}\alpha^2 \frac{\int |\nabla
h|^2}{\int |\nabla Q|^2} \right) \left(1+\alpha \frac{\int Q
h_1}{\int Q^2}+\frac 12\alpha^2\frac{\int |h|^2}{\int Q^2}
-\frac 12\alpha^2\frac{(\int Q h_1)^2}{\left(\int Q^2\right)^2}\right)\\
-\left(1+\alpha \frac{\int 4 Q^3 h_1}{\int Q^4} +\alpha^2\frac{\int
6 Q^2 h_1^2+2Q^2 h_2^2}{\int Q^4}\right)+O\big(\alpha^3\big).
\end{multline*}
Since $I(Q)=0$ and $I(Q+\alpha h)\geq 0$ for all real $\alpha$, the
linear term in $\alpha$ in the previous estimate is zero, and the
quadratic term is nonnegative. Applying $\int|\nabla Q|^2=3\int Q^2$
and $\int Q^4=4\int Q^2$, we get
\begin{equation*}
\frac{3\int |\nabla h|^2}{2\int |\nabla Q|^2}+\frac{\int |
h|^2}{2\int Q^2}-\frac{6\int Q^2 h_1^2}{\int Q^4}-\frac{2\int Q^2
h_2^2}{\int Q^4}-\frac{\left(\int Q h_1\right)^2}{2\left(\int
Q^2\right)^2}= \frac{\Phi(h)}{\int Q^2}-\frac{\left(\int Q
h_1\right)^2}{2\left(\int Q^2\right)^2}\geq 0,
\end{equation*}
which implies \eqref{positivity}.
\medskip

\noindent\emph{Step 2. Coercivity.} We show that if $h$ fulfils
\eqref{ortho1} and \eqref{ortho2}, then for some $c_*$,
\begin{equation}
\label{coercivity2}
\Phi(h)\geq c_{*} \|h\|_{H^1}^2.
\end{equation}
Note that $\Phi(h)=\Phi_1(h_1)+\Phi_2(h_2)$, where
\begin{align*}
\Phi_1(h_1)&:=\frac{1}{2} \int |\nabla h_1|^2+\frac{1}{2} \int h_1^2
- \frac{3}{2} \int Q^2 h_1^2=\frac 12\int (L_+h_1)\,h_1,\\
\Phi_2(h_2)&:=\frac{1}{2} \int |\nabla h_2|^2+\frac{1}{2} \int h_2^2
- \frac{1}{2} \int Q^2 h_2^2=\frac 12\int (L_-h_2)\,h_2.
\end{align*}
By Step 1, $L_+$ is nonnegative on $\{\Delta Q\}^{\bot}$ and $L_-$
is nonnegative. We will deduce \eqref{coercivity2} from Remark
\ref{R:null.space} and a classical argument (see the proof of
Proposition 2.9 in \cite{We85}).

We claim that under assumptions \eqref{ortho1} and \eqref{ortho2},
there exists $c>0$ such that $\Phi_1(h_1)\geq c\,\|h_1\|_{H^1}^2$.
For this we first show that there exists $c_1>0$ such that
\begin{equation}
\label{coercivity3}
\int \partial_{x_1}Q\,h_1= \int \partial_{x_2}Q\,h_1= \int \partial_{x_3}Q\,h_1
=\int \Delta Q\, h_1=0\Longrightarrow  \Phi_1(h_1)\geq c_1\|h_1\|_{2}^2.
\end{equation}
Following the proof of \cite[Proposition 2.9]{We85}, assume that
\eqref{coercivity3} does not hold. Then there exists a sequence of
real-valued $H^1$-functions $\{f_n\}_n$ such that
\begin{gather}
 \label{absurd.coer}
\lim_{n\rightarrow +\infty} \Phi_1(f_n)=0,\quad \|f_n\|_{2}=1 \quad \text{and}\\
 \label{absurd.ortho}
\int \Delta Q f_n=\int (\partial_{x_1}Q) f_n=\int (\partial_{x_2}Q)
f_n=\int (\partial_{x_3}Q) f_n=0.
\end{gather}
In particular,
\begin{equation}
 \label{absurd.coer2}
\frac 12 \int |\nabla f_n|^2 =-\frac 12\int |f_n|^2 +\frac 32 \int
Q^2 f_n^2 + o(1) =-\frac 12+\frac 32 \int Q^2 f_n^2+o(1).
\end{equation}
Thus, $\int |\nabla f_n |^2 \leq C \, \|f_n\|_2=C$, and hence,
$\{f_n\}_n$ is bounded in $H^1$. Extracting, if necessary, a
subsequence from $\{f_n\}_n$, we get that there exists $f_*\in H^1$
such that
\begin{equation}
 \label{weakCV}
f_n \underset{n\rightarrow+\infty}{\relbar\!\rightharpoonup}
f_*\text{ weakly in }H^1.
\end{equation}
Since $Q$ is decreasing at infinity, we have
\begin{equation}
 \label{stronglimfn}
\frac{3}{2}\int Q^2 f_n^2\underset{n\rightarrow+\infty}{\longrightarrow}
\frac{3}{2}\int Q^2 f_*^2.
\end{equation}
By \eqref{absurd.coer2}, it follows that
\begin{equation}
\label{f*nonzero}
\int Q^2 f_*^2\geq \frac 13,
\end{equation}
and, in particular, $f_*\neq 0$. Furthermore, $\limsup_n
\|f_n\|_{H^1}\leq \|f_*\|_{H^1}$, and thus, $\limsup_n
\Phi_1(f_*)\leq \limsup_n \Phi_1(f_n)$. By \eqref{absurd.coer},
$\limsup_n \Phi_1(f_*)\leq 0$. By the weak convergence, $\int \Delta
Q f_*=\int (\partial_{x_1}Q) f_*=\int (\partial_{x_2}Q) f_*=\int
(\partial_{x_3}Q) f_*=0$. In particular, by Step 1, $\Phi_1(f_*)\geq
0$.  Therefore,
\begin{equation}
 \label{Phi1f*}
\frac 12 \int (L_+f_*) \,f_*=\Phi_1(f_*) =0,
\end{equation}
and $f_*$ is the solution to the following minimization problem
\begin{multline*}
0=\frac{\int (L_+f_*)\,f_*}{\|f_*\|_2}=\min_{f\in E\setminus \{0\}}
\frac{\int (L_+f)\, f}{\|f\|_2},\quad f_*\in E, ~~\text{where}\\
 E:=\left\{f\in H^1,\;\int \Delta Q f=\int (\partial_{x_1}Q) f
 =\int (\partial_{x_2}Q) f=\int (\partial_{x_3}Q)f=0\right\}.
\end{multline*}
Hence, for some Lagrange multipliers $\lambda_0,\ldots,\lambda_3$,
we can write
\begin{equation}
 \label{L+f*}
L_+ f_* =\lambda_0 \Delta Q +\lambda_1 \partial_{x_1}
Q+\lambda_2\partial_{x_2}Q +\lambda_3\partial_{x_3}Q.
\end{equation}
By the symmetry of $Q$,
$$
\int (\partial_{x_i}Q)\partial_{x_j}Q=0,\;i\neq j, \quad \int
(\Delta Q)\partial_{x_i}Q=0.
$$
By Remark \ref{R:null.space}, we have
$L_+\partial_{x_i}Q=0$ for $i\in\{1,2,3\}$, which from \eqref{L+f*}
gives
$$
0=-\int f_*L_+(\partial_{x_i}Q)=\int (L_+ f_*)\partial_{x_i}Q =
\lambda_i\int\left|\partial_{x_i}Q\right|^2,
$$
showing that $\lambda_1=\lambda_2=\lambda_3=0$. Hence,
\begin{equation}
\label{L1f*}
L_+f_*=\lambda_0 \Delta Q=\lambda_0\left(-Q^3+Q\right).
\end{equation}
Denote $\widetilde{Q}=Q+x\cdot \nabla Q$. Let us show that
\begin{equation}
 \label{E:L+f*}
L_+\left(\frac{\lambda_0}{2}(Q-\wtQ)\right)=\lambda_0(-Q^3+Q)=L_+ f_*.
\end{equation}
Indeed, $L_+ Q=-2Q^3$. Furthermore, if $Q_{\lambda}(x)=\lambda
Q(\lambda x)$, then $ \wtQ:=\frac{\partial}{\partial
\lambda}\left(Q_{\lambda} \right)_{\restriction \lambda =1}$.
Differentiating the equality $-\lambda^2Q_{\lambda}+\Delta
Q_{\lambda} +Q_{\lambda}^3 = 0$ with respect to $\lambda$ at
$\lambda=1$, we obtain $L_+ \wtQ=-2Q$, which produces
\eqref{E:L+f*}.

By Remark \ref{R:null.space} and \eqref{E:L+f*}, there
exist $\mu_1$, $\mu_2$, $\mu_3$ such that
$$
f_* = \frac{\lambda_0}{2}(Q -\wtQ) + \mu_1 \, \partial_{x_1} Q +
\mu_2 \, \partial_{x_2}Q + \mu_3 \, \partial_{x_3}Q.
$$
Next, note that $\int \wtQ\partial_{x_j}Q=0$ (indeed, $\int
Q\partial_{x_j}Q=0$ by integration by parts, and $\int
x_i\partial_{x_i}Q\partial_{x_j}Q=0$ by the symmetry of $Q$). Using
that $\int f_*\partial_{x_j}Q=0$, we get $ \mu_1=\mu_2=\mu_3=0$.
Hence,
\begin{equation}
\label{f*}
f_*=\frac{\lambda_0}{2}(Q-\wtQ)=-\frac{\lambda_0}{2} x\cdot\nabla Q.
\end{equation}
By straightforward calculation, \eqref{L1f*} and \eqref{f*}, we
obtain
$$
\Phi_1(f_*) = \frac{1}{2}\int (L_+ f_*) f_* = -\frac{\lambda_0^2}{4}
\int \Delta Q (x\cdot \nabla Q) = -\frac{\lambda_0^2}{8} \int
|\nabla Q|^2.
$$
By \eqref{Phi1f*}, $\lambda_0=0$, and therefore, $f_*=0$, which
contradicts \eqref{f*nonzero} and concludes the proof of
\eqref{coercivity3}.

By the explicit expression of $\Phi_1$, we have that for $\eps>0$
small enough, $\eps\Phi_1(h_1)\geq \frac{\eps}{2} \int |\nabla
h_1|^2-\frac{c_1}{2}\int h_1^2$ for any $h_1\in H^1$, where $c_1$ is
the constant in \eqref{coercivity3}. Adding to \eqref{coercivity3},
we get that for some constant $c>0$,
$$
\int (\partial_{x_1}Q)h_1=\int (\partial_{x_2}Q)h_1 = \int
(\partial_{x_3}Q) h_1 = \int \Delta Qh_1=0  \Longrightarrow
\Phi_1(h_1)\geq c\|h_1\|_{H^1}^2.
$$
To conclude the proof of Proposition \ref{P:Coercivity}, it
remains to show that for some constant $c>0$,
$$
\int Q h_2=0 \Longrightarrow \Phi_2(h_2)\geq c\|h_2\|_{H^1}^2.
$$
The proof is similar to the previous and we omit it.

\subsection{Coercivity of $\Phi$ on $G_{\bot}'$}

We first show
\begin{equation}
\label{E:G'>0}
\forall h\in G_{\bot}'\setminus \{0\},\quad \Phi(h)>0.
\end{equation}
If not, there exists $\tilde{h}\in
H^1\setminus \{0\}$ such that
\begin{equation}
\label{AbsurdOrtho} \int (\partial_{x_j}Q)
\tth_1=0,\;j=1,2,3,\quad\int Q \tth_2 = \int \YYY_1\tth_2=\int
\YYY_2 \tth_1=0\quad\text{and}\quad \Phi\big(\tth\big)\leq 0.
\end{equation}
Recall that by Remark \ref{R:null.space},
$$
\forall \, h\in H^1, \quad B(\partial_{x_1}Q,h) =
B(\partial_{x_2}Q,h) = B(\partial_{x_3}Q,h)=B(iQ,h)=0.
$$
Furthermore, by \eqref{AbsurdOrtho},
\begin{equation}
\label{Y+th}
B(\YYY_+, \tth) = \frac12 \, \int (L_+\YYY_1) \, \tth_1 + \frac12 \,
\int (L_-\YYY_2) \, \tth_2 = \frac12\, e_0 \int \YYY_2 \,
\tth_1-\frac12\,  e_0 \int \YYY_1\, \tth_2=0,
\end{equation}
so we have that
$\partial_{x_j}Q$, $j=1,2,3$, $iQ$, $\YYY_+$ and $\tth$ are
orthogonal in the bilinear symmetric form $B$. Noting that
$\Phi(iQ)=\Phi(\partial_{x_j}Q)=\Phi(\YYY_+)=0$ and that
$\Phi(\tth)\leq 0$, we get
$$
\forall h \in E:= \vect\{\partial_{x_1} Q,\partial_{x_2}Q,
\partial_{x_3}Q,iQ,\YYY_+,\tth\}, \quad \Phi(h)\leq 0.
$$
We claim that $\dim_{\RR} E=6$.

Assume that for some real numbers $\alpha_j$, $\beta$, $\gamma$,
$\delta$, we have
\begin{equation}
\label{DependentVectors} \sum_{j=1}^3 \alpha_j \partial_{x_j}Q+\beta
iQ +\gamma \YYY_+ +\delta \tth=0.
\end{equation}
By Remark \ref{R:L-positive}, $B(\YYY_+,{\YYY_-}) = -e_0 \,
(L_-\YYY_2,\YYY_2) \neq 0$. Furthermore, the same computation as in
\eqref{Y+th} shows that  $B(\tth,\YYY_-)=0$. From
\eqref{DependentVectors} we get that $\gamma
\,B(\YYY_+,{\YYY}_-)=0$, which implies $\gamma=0$. Since
$\partial_{x_j}Q$, $iQ$ and $\tth$ are orthogonal in $L^2$, we also
get that $\alpha_j=\beta=\delta=0$. Thus, $\dim E=6$.

We know that $\Phi$ is definite positive on $G_{\bot}$ (a subspace
of codimension $5$ of $H^1$), hence, cannot be non-positive on $E$,
yielding a contradiction. The proof of \eqref{E:G'>0} is complete.

It remains to show that if $h\in G_{\bot}'$, $\Phi(h)\geq
c\|h\|^2_{H^1}$. Let us sketch the proof. As before, it is
sufficient to show
\begin{equation}
\label{E:coercivity_L2_G'}
\exists c>0,\;\forall h\in G_{\bot}',\quad \Phi(h)\geq c\|h\|_{2}^2.
\end{equation}
If not, there exists a sequence $h_n\in G_{\bot}'$ such that
\begin{equation}
\label{E:absurd_L2_G'}
 \lim_{n\rightarrow +\infty} \Phi(h_n)=0,\text{ and } \forall n,\;\|h_n\|_{2}^2=1.
\end{equation}
Extracting a subsequence from $(h_n)$ if necessary, we may assume
$h_n \rightharpoonup h^*$ weakly in $H^1$. The weak convergence of
$h_n$ implies $h^*\in G_{\bot}'$. By \eqref{E:absurd_L2_G'} it is
easy to check that $h^*\neq 0$ and $\Phi(h^*)=0$, which contradicts
\eqref{E:G'>0}, showing as announced \eqref{E:coercivity_L2_G'}.
\qed

\section{Proof of a Cauchy-Schwarz type inequality}
\label{A:claimval} Let us prove Claim \ref{C:val}. Let
$\dd(f)=\int |\nabla Q|^2 -\int |\nabla f|^2$ (so that
$\delta(f)=|\dd(f)|$) and  $\lambda\in \RR$. Then
$\left\|e^{i\lambda \varphi}f\right\|_2 =\| Q \|_2$. By the
Gagliardo-Nirenberg inequality
$$
\left\| \nabla \big(e^{i\lambda \varphi} f \big) \right\|_2^3
\|Q\|_4^4 \geq \|\nabla Q\|_2^3 \, \|f\|_4^4.
$$
Raising the previous inequality to the power $2/3$, and expanding
$\left\|\nabla \left(e^{i\lambda \varphi}f\right)\right\|_2^2$, we
get
$$
\lambda^2 \int |\nabla \varphi|^2 |\nabla f|^2+ 2\lambda \im \int
(\nabla \varphi \cdot \nabla f) \, \fbar + \int |\nabla f|^2 -
\|f\|_4^{\frac 83} \frac{\|\nabla Q\|_2^{2}}{\|Q\|_4^{\frac 83}}
\geq 0.
$$
Using elementary properties of quadratic inequalities (in
$\lambda$), we obtain
$$
\left|\im \int (\nabla \varphi \cdot \nabla f) \, \fbar \right|^2 \leq
\left(\int |\nabla \varphi|^2 |\nabla f|^2\right)\left(\int |\nabla
f|^2 -\|f\|_4^{\frac 83}\frac{\|\nabla Q\|_2^{2}}{\|Q\|_4^{\frac
83}} \right).
$$
We have $\|\nabla f\|^2_2 =\|\nabla Q\|^2_2-\dd(f)$ and, by
assumption \eqref{hyp.val}, we have $ \| f \|_4^4 = \| Q \|_4^4 - 2
\, \dd(f)$, so that
\begin{align*}
\int |\nabla f|^2_2 -\|f\|_4^{\frac 83}\frac{ \|\nabla
Q\|_2^{2}}{\|Q\|_4^{\frac 83}} & =\|\nabla
Q\|_2^2-\dd(f)-\left(\|Q\|_4^4-2\dd(f)\right)^{\frac 23}
\frac{\|\nabla Q\|_2^{2}}{\|Q\|_4^{\frac 83}}\\
& =\|\nabla Q\|_2^2-\dd(f)-\|\nabla Q\|_2^2+\frac 43
\frac{\dd(f)}{\|Q\|_4^4}\|\nabla Q\|_2^2+O\left(|\dd^2(f)|\right).
\end{align*}
Recalling that $\|Q\|_4^4=\frac 43\|\nabla Q\|_2^2$, we obtain
$$ \int |\nabla f|^2_2 -\|f\|_4^{\frac 83}\frac{ \|\nabla
Q\|_2^{2}}{\|Q\|_4^{\frac 83}}=O\left(|\dd^2(f)|\right),$$
concluding the proof of the claim.
\qed


\end{document}